\documentstyle[amssymb, 11pt]{article}

\catcode`@=11
\def\relaxnext@{\let\next\relax}
\def\noaccents@{\def\accentfam@{0}}
\font\teneuf=eufm10
\font\seveneuf=eufm7
\font\fiveeuf=eufm5
\newfam\euffam
\textfont\euffam=\teneuf
\scriptfont\euffam=\seveneuf
\scriptscriptfont\euffam=\fiveeuf
\def\frak{\relaxnext@\ifmmode\let\next\frak@\else
 \def\next{\Err@{Use \string\frak\space only in math mode}}\fi\next}
\def\goth{\relaxnext@\ifmmode\let\next\frak@\else
 \def\next{\Err@{Use \string\goth\space only in math mode}}\fi\next}
\def\frak@#1{{\frak@@{#1}}}
\def\frak@@#1{\noaccents@\fam\euffam#1}
\catcode`@=12

\newcommand{\reals}{{\Bbb R}}

\newcommand{\complex}{{\Bbb C}}
\newcommand{\T}{{\Bbb T}}

\newcommand{\fg}{{\frak g}}
\newcommand{\fk}{{\frak k}}
\newcommand{\ft}{{\frak t}}

\newcommand{\s}{{\star}}

\newcommand{\CA}{{\cal A}}

\newcommand{\CL}{{\cal L}}

\newcommand{\CH}{{\cal H}}
\newcommand{\CX}{{\cal X}}

\newcommand{\DIF}{{\Omega^1}}
\newcommand{\CXP}{{\cal X}^{\star}(P)}

\newcommand{\OP}{\Omega^{\star}(P)}
\newcommand{\HP}{H^{\star}_{\pi}(P)}

\newcommand{\Cas}{\mbox{\em Cas}}

\newcommand{\ta}{\tilde{a}}

\newtheorem{thm}{Theorem}[section]
\newtheorem{prop}[thm]{Proposition}
\newtheorem{lemma}[thm]{Lemma}
\newtheorem{cor}[thm]{Corollary}
\newtheorem{dfn}[thm]{Definition}
\newtheorem{rmk}[thm]{Remark}
\newtheorem{exam}[thm]{Example}

\newcommand{\qed}{\begin{flushright} $\Box$\ \ 
                 \end{flushright}}

\newcommand{\proof}{\noindent{\em Proof}.~}
\newcommand{\pproof}{\noindent{\em Proof}~}

\title{{\bf Equivariant Poisson Cohomology and a 
Spectral Sequence Associated with a Moment Map}}
\author{{\bf Viktor L.  Ginzburg}\thanks{The author was supported in 
part by the National Science Foundation}
\thanks{Current address:
Department of Mathematics, 
University of California at Santa Cruz,
Santa Cruz, CA 95064, USA;
ginzburg@cats.ucsc.edu}\\
\\
Department of Mathematics\\
University of California at Berkeley\\
Berkeley, CA 94720\\
ginzburg@math.berkeley.edu}
\date{November, 1996}

\begin{document}

\maketitle

\begin{abstract}
We introduce and study a new spectral sequence 
associated  with a Poisson group action on a Poisson manifold and 
an equivariant momentum mapping.
This spectral sequence is a Poisson analog of the Leray spectral sequence 
of a fibration.
The spectral sequence converges to the Poisson cohomology of
the manifold and has the $E_2$-term equal to the tensor product of 
the cohomology of the Lie algebra and the equivariant Poisson 
cohomology of the manifold. The latter is defined as
the equivariant cohomology of the multi-vector fields made into a 
$G$-differential complex by means of the momentum mapping. An extensive
introduction to equivariant cohomology of $G$-differential complexes is
given including some new results and a number of examples and applications
are considered.

\end{abstract}

\tableofcontents{}

\section{Introduction}

The main result of the present paper is the existence of
a spectral sequence associated  with an equivariant momentum mapping of a 
Poisson manifold and converging to the Poisson cohomology
of the manifold. When the manifold is symplectic, the
spectral sequence is isomorphic to that in the de Rham
cohomology for the dual data: the group action rather
than the momentum mapping. Namely, let $P$ be a symplectic 
manifold acted on by a compact Poisson Lie group $G$. Then
the spectral sequence in question is isomorphic
to the Leray spectral sequence of the principal $G$-bundle
$P\times EG\to (P\times EG)/G$. If the
momentum mapping is a submersion onto an open subset, and so
the action is locally free, one can take the projection
$P\to P/G$ as well. 

To extend the spectral sequence to Poisson manifolds, one thus
needs to dualize the notion of equivariant de Rham cohomology,
which leads to what we call equivariant Poisson cohomology. The 
latter arises as a particular case of the equivariant cohomology
of a $G$-differential complex, for an equivariant momentum
mapping makes the complex of multi-vector fields into one.

Equivariant de Rham cohomology has proven itself to be a very convenient
tool in the study of Hamiltonian group actions because it renders
a powerful topological machinery to one's service. However, no analog of
such a cohomology has been known for Poisson manifolds. 
The introduction of equivariant Poisson cohomology fills in the
gap. The cohomology readily lends itself to a variety of applications
(e.g., the spectral sequence and calculations of Poisson cohomology)
in a way similar to de Rham cohomology, but with algebraic topology 
replaced by some algebra. Certain aspects of the new techniques
have remained unexplored in this paper. For example, one may see that
the equivariant Poisson cohomology spaces behave well under symplectic
reduction analogously to equivariant de Rham cohomology. However, no 
rigorous results in this direction have been proved yet. It is also worth 
noticing that, although the formal algebraic machinery has been developed 
long time ago, Poisson equivariant 
cohomology is still one of very few examples, if not the only example, where
$G$-differential complexes essentially different from the de Rham complex 
are employed

The spectral sequence introduced in the paper can be generalized to virtually
any Poisson map. The generalization is particularly transparent for 
Poisson submersions. This gives an affirmative
answer to a question asked by Alan Weinstein whether there exists a 
spectral sequence associated with a ``Poisson fibration''.
In general, the $E_1$ and $E_2$-terms tend to have a very complicated
structure. However, for a momentum mapping, the spectral
sequence is relatively easy to analyze and, as illustrated by
a number of examples, use to calculate the Poisson cohomology, to which 
it converges. 

Before the spectral sequence can be defined and the
$E_2$ calculated,  a rather cumbersome machinery has to be 
developed. Some work to this end has already been carried out in 
\cite{gi:moment}. The exposition in the present paper naturally evolves 
around the 
study of equivariant cohomology of $G$-differential complexes and its
Poisson counterpart. The paper is organized as follows. 

We begin with an extensive introduction, enlivened with examples,
to  equivariant cohomology of $G$-differential complexes. This approach,
developed originally by H. Cartan \cite{cartan1}, is based on the 
axiomatization within an algebraic framework of the key features of
differential forms needed to define
the equivariant cohomology of a manifold with real coefficients.
Thus a $G$-differential complex is, by definition, a complex $A^\s$
equipped with a $G$-action and contractions with the elements 
of $\fg$ satisfying the natural compatibility conditions. A fundamental
example is the de Rham complex $A^\s=\Omega^\s(M)$, where $M$ is
a manifold acted on by a group $G$. The equivariant cohomology
$H^\s_G(A^\s)$ is then just the cohomology of the complex
$(A^\s\otimes S^\s\fg^*)^G$, which for the de Rham complex
gives $H^\s_G(M,\reals)$, provided that $G$ is compact.

In the opening subsections of Section \ref{sect:poisson-general}, we 
give all the necessary definitions and state first basic theorems including
the result that for ``locally free $G$-differential 
algebras'' the equivariant cohomology coincides with the basic cohomology.
An alternative source for this material is Appendix B of 
\cite{gls}  or the second part of \cite{DKV}. 

Then we turn to one of our main examples of $G$-differential complex -- 
the Chevalley--Eilenberg complex $C^\s(\fg;V)$, where $V$ is a Fr\'{e}chet
$G$-module. This is a complex over a $G$-differential locally free algebra
$C^\s(\fg)$, and so $H^\s_K(\fg;V)=H^\s(\fg,\fk;V)$ for any subgroup 
$K\subset G$. In Section \ref{section:lie-alg}, we recall that
$H^\s(\fg,\fk;V)=H^\s(\fg,\fk)\otimes V^G$, when $G$ is compact,
which is rather easy to prove using
the so-called van Est spectral sequence (cf., \cite{gi:moment}). A
proof that does not rely on the spectral sequence and gives a little 
bit stronger parametric version of the theorem is provided in the appendix.
The results of this section play a key role in the explicit calculation
of (equivariant) Poisson cohomology of certain Poisson manifolds, e.g., 
$G^*$, and of the $E_2$-term of the spectral sequence. Conceptually, this
role is not entirely dissimilar to that of Hodge theory in analysis on 
manifolds.

In Section \ref{subsec:further} we study further properties of equivariant
cohomology. For example,
we analyze the equivariant cohomology in low degrees and introduce
a spectral sequence analogous to the standard spectral sequence for 
equivariant cohomology, and thus having $E_2=H^\s(A^\s)\otimes (S^\s\fg^*)^G$
and converging to $H^\s_G(A^\s)$.

Finally, in Section \ref{sec:spec-general} we call into action 
the major figure of this paper -- a spectral sequence relating
the cohomology of $A^\s$, the cohomology of $\fg$, and the equivariant
cohomology $A^\s$. As we have mentioned, the spectral sequence is a 
generalization of the Leray spectral sequence of 
$M\times EG\to (E\times EG)/G$, where $M$ is a space acted on by $G$.
It has $E_2=H^\s_G(A^\s)\otimes H^\s(\fg)$ and converges to the
ordinary cohomology $H^\s(A^\s)$.

Section \ref{sec:dfn} is a review of the results from Poisson geometry
used in the subsequent sections. Along
with the standard material, which can be found elsewhere (see, e.g., 
\cite{va:book} and bibliography therein), we prove some new results
and outline a technique developed in \cite{gi:moment} to be applied
in the present paper.

Section \ref{sec:Poisson-cohomology} is entirely devoted to equivariant
Poisson cohomology. The data sufficient to define this cohomology,
$H^\s_{\pi,G}(P)$, are the Poisson manifold $P$ and 
an equivariant momentum mapping $\mu\colon P\to G^*$, and so a Poisson
action of a Poisson Lie group $G$. We start with a discussion of the 
general properties of the cohomology and carry out some calculations in 
low degrees. The ``locally free'' case where $\mu$ is a submersion onto
an open subset on $G^*$ is a key to understanding the geometrical meaning
of the cohomology. In Section \ref{subsec:tangent}, we show that in this
case, $H^\s_{\pi,G}(P)$ can be identified with the cohomology of 
$G$-invariant multi-vector fields tangent to the $\mu$-fibers.
(It is instructive to contrast this result with the identification
of the equivariant and basic de Rham cohomology for locally free actions.) 
The equivariant cohomology of $G^*$ with respect to a
subgroup $K\subset G$ is found in Section \ref{subsec:calcul}:
$H^\s_{\pi,K}(G^*)=H^\s(\fg,\fk)\otimes C^\infty(G^*)^G$, and, in 
particular, $H^\s_{\pi,G}(G^*)= C^\infty(G^*)^G$. Little proving needs
to be done throughout Section \ref{sec:Poisson-cohomology}, for the
results proved previously for general $G$-differential complexes readily 
apply here.

Finally, in Section \ref{sec:spec-mom} we define the spectral sequence
associated with a momentum mapping, calculate its $E_2$-term and analyze
a variety of examples. 

Throughout the paper  the cohomology of a manifold are taken with
real coefficients unless specified otherwise. All Lie groups are
assumed to be connected.

{\bf Acknowledgments.} I would like to thank 
 Mich\'{e}le Vergne, Alan Weinstein, and Gregg Zuckerman 
for a number of useful remarks. I am also grateful to
the Isaac Newton Institute for its kind hospitality
and providing a stimulating environment during the 1994 program on symplectic
geometry when a part of the manuscript was prepared.

\section{Equivariant cohomology}
\label{sect:poisson-general}

In this section we study the equivariant cohomology in the general 
algebraic context of $G$-differential complexes, introduced originally in 
by H. Cartan in \cite{cartan1}. 
The material reviewed in the first two subsections is fairly standard
and the proofs are omitted. The reader interested in more details
should consult, e.g., \cite{DKV} or Appendix  B of \cite{gls}.
(Many other relevant references are spread throughout the text.)
Two concluding subsections contain some newer results crucial for 
what follows and their detailed proofs are provided.

\subsection{Equivariant cohomology of $G$-differential complexes}
\label{subsec:G-differential}

Let $(A^{\s}, d)$ be a complex of Fr\'{e}chet spaces with the differential
$d$ of degree one. For the sake of simplicity we assume throughout this 
section that the grading on $A^\s$ is positive: $A^n=0$ when $n<0$. 
Assume that a connected Lie group $G$ acts 
smoothly on $A^{\s}$ and that the action commutes with $d$. Denote by 
$L_{\xi}$ the infinitesimal action on $A^{\s}$ of an element $\xi$ of the 
Lie algebra $\fg$ of $G$. The complex $A^{\s}$ is said to be
a {\em $G$-differential complex} if for any $\xi\in\fg$ it is equipped with a
continuous linear mapping $i_{\xi}\colon A^{\s}\to A^{\s}$ of degree -1 
such that
the following conditions hold for all $\xi$ and $\zeta$ in $\fg$ and  $g\in G$:
\begin{enumerate}
\item[(i)] $i_{\xi}i_{\zeta}+i_{\zeta}i_{\xi}=0$~;

\item[(ii)] $gi_{\xi}g^{-1}=i_{ad_g\xi}$~;

\item[(iii)] $L_{\xi}=di_{\xi}+i_{\xi}d$ (Cartan's identity).
\end{enumerate}
If instead of a  $G$-action, $A^\s$ carries only an infinitesimal 
$\fg$-action, as above, commuting with the differential, $A^{\s}$ is 
said to be a
{\em $\fg$-differential complex}, provided that we have ``contractions''
$i_{\xi}$, (i) and (iii) hold, and (ii) is replaced by its
infinitesimal version:
\begin{enumerate}
\item[(ii')] $i_{[\xi,\zeta]}=L_{\xi}i_{\zeta}-i_{\zeta}L_{\xi}$~.
\end{enumerate}
Let $A^{\s}$ be a graded-commutative differential algebra, i.e.,
$ab=(-1)^{\deg a\deg b}ba$, and let
the graded Leibniz identity hold for $d$ and the contractions, i.e.,
$i_{\xi}(ab)=(i_{\xi}a)b+(-1)^{\deg a}a i_{\xi}b$.
Then we call $A^{\s}$ a $G-$ or {\em $\fg$-differential algebra}.

\begin{rmk}
{\em
Our definition of a $\fg$-differential complex is redundant in the sense
that the action of $\fg$ on $A^{\s}$ can be recovered from $d$ and
the ``contractions''. To be more precise, assume that $A^{\s}$ is a
complex and $i_{\xi}$ are given so that (i) and (ii') hold. Define the 
infinitesimal action $L_{\xi}$ via Cartan's identity (iii). The
commutation relation $L_{[\xi,\zeta]}=L_{\xi}L_{\zeta}-L_{\zeta}L_{\xi}$
is then a consequence of (i) and (ii'). Furthermore, if $A^\s$ is a
$G$- or $\fg$-differential algebra, the multiplication is 
automatically $G$- or $\fg$-invariant .

Note also that a $\fg$-differential complex may fail to integrate to
a $G$-differential complex even if $G$ is compact.
}
\end{rmk}

Consider the tensor product $A^{\s}\otimes S^{\s}\fg^*$
with the grading $\deg(A^n\otimes S^m\fg^*)=n+2m$. This
tensor product 
can be viewed as the space of $A^{\s}$-valued polynomial
functions on $\fg$. Thus for $\alpha\in A^{\s}\otimes S^{\s}\fg^*$
and $\xi\in \fg$ we denote by $\alpha(\xi)\in A^{\s}$ the value of
such a polynomial at $\xi$. Let us define a homomorphism 
$d_G\colon A^{\s}\otimes S^{\s}\fg^*\to A^{\s}\otimes S^{\s}\fg^*$
of  degree one as follows:
\begin{equation}
\label{eq:equi-diff}
(d_G\alpha)(\xi)=d(\alpha(\xi))-i_{\xi}\alpha(\xi)
\enspace . 
\end{equation}
A routine calculation shows that $d^2_G\alpha(\xi)=-L_{\xi}\alpha(\xi)$,
and so $d_G$ is not a differential on the entire tensor product.
Note however that $A^{\s}\otimes S^{\s}\fg^*$ is a $\fg$-module with the 
diagonal action  and denote by $A^{\s}_{\fg}$  its subspace 
$(A^{\s}\otimes S^{\s}\fg^*)^{\fg}$ of $\fg$-invariant
elements. Clearly, $A^{\s}_{\fg}$ is closed under $d_G$ and
$d_G^2=0$ on $A^{\s}_{\fg}$. 
\begin{dfn}
\label{def:equi-general}
{\em
The cohomology of the complex $(A^{\s}_{\fg}, d_G)$,
denoted $H^{\s}_G(A^\s)$, is called  the {\em
equivariant cohomology} of the $\fg$-differential complex $A^{\s}$.
The complex $(A^{\s}_{\fg}, d_G)$ is said to be the {\em Cartan model}
for equivariant cohomology.
}
\end{dfn}

Let us emphasize that in what follows we use the notation $H^\s_G(A^\s)$
even when $A^\s$ carries only an infinitesimal $G$-action, i.e.,
when $A^\s$ is a $\fg$-module.

Recall that
since $d_G$ is linear over $(S^{\s}\fg^*)^{\fg}$, the equivariant
cohomology is a graded $(S^{\s}\fg^*)^{\fg}$-module. If $A^\s$
is a $\fg$-differential algebra,  the multiplication descends
to $H^\s_G(A)$ making $H^\s_G(A^\s)$ into an $(S^{\s}\fg^*)^{\fg}$-algebra.

\begin{exam}
\label{exam:de-rham}
{\em
The standard example of equivariant cohomology is the
equivariant de Rham cohomology. Let $M$ be a manifold acted on by $G$.
We take  the de Rham complex of $M$ with the natural action of $G$ and 
natural contractions as the $G$-differential algebra 
$(A^\s,d)=(\Omega^\s(M),d)$.
Then the cohomology $H^\s_G(A^\s)$ is called the equivariant cohomology
of $M$ with real coefficients. To be more precise, let
$H_G^\s(M,\reals)$ be the ordinary (real) cohomology of the
quotient $(M\times EG)/G$ where $EG$ is the classifying space of $G$. Then
$H^\s_G(M,\reals)=H^\s_G(A)$ when $G$ is compact (the equivariant de Rham
theorem).

The $\fg$-differential algebra construction goes through even 
if $G$ is not compact or if $M$ is equipped 
only with an infinitesimal $G$-action. In these cases,
the equivariant cohomology of $A$ will also be called the equivariant 
cohomology of $M$ even though the topological construction may then lead
to an entirely different result \cite{at-bo:moment_map}.

Likewise, let $E\to M$ be a vector bundle (with a Fr\'{e}chet fiber)
endowed with a flat connection. We use horizontal lifts of vector vector fields
to obtain a $G$-action on $E$. Then the complex $A^\s=\Omega^\s(M,E)$ of
$E$-valued differential forms turns naturally into a $G$-differential
complex, and we have the equivariant cohomology $H^\s_G(M;E)$. 
When $G$ is compact, this space is just the standard equivariant
cohomology of $M$ with coefficients in the sheaf (a local system) of
flat sections of $E$. Note that now $A^\s$ is not a $G$-differential 
algebra. However, it will be such when the fibers of $E$ are algebras
and multiplication is invariant of holonomy.

The comparison of algebraic and topological constructions
of equivariant cohomology will be continued in Example 
\ref{exam:de-rham-loc-free}.
}
\end{exam}

\subsection{The Weil complex and locally free $G$-differential complexes}
\label{sec:Weil}

In this section we outline a slightly different, though equivalent, approach  
due to Atiyah and Bott \cite{at-bo:moment_map} to the
definition of equivariant cohomology.
This approach relying on the usage of the Weil algebra of $\fg$ appears to be 
considerably more convenient for our present purposes.
However, it should be mentioned that adapting the methods of \cite{cartan2}
one can prove Theorem \ref{thm:basic-complex} and Corollary 
\ref{thm:basic-algebra} without ever leaving the realm of complexes 
$A^{\s}_{\fg}$.

Let $W^\s(\fg)$ be the Weil algebra of $\fg$.
Recall that, as a graded commutative algebra, 
$W^\s(\fg)=\wedge^\s\fg^*\otimes S^\s\fg^*$ with the symmetric part
given the even grading. The elements of $W^\s(\fg)$ can be thought of
either as symmetric (even) functions on $\fg$ with values in
$\wedge^\s\fg^*$ or odd (skew-symmetric) functions on $\fg$ with
values in $S^\s\fg^*$. The second interpretation leads to the 
$S^\s\fg^*$-linear contractions 
$i_\xi\colon W^\s(\fg)\to W^\s(\fg)$ for $\xi\in\fg$.
To be more precise, $i_\xi(\varphi\otimes f)=(i_\xi\varphi)\otimes f$, where 
$\varphi\in\wedge^\s\fg^*$ and $f\in S^\s\fg^*$. 

Furthermore, $W^\s(\fg)$ carries a differential $d_W=d_L-\delta$, 
Here $d_L$ is the 
Chevalley--Eilenberg differential of $\fg$ with coefficients in 
the $\fg$-module $E=S^\s\fg^*$ (see (\ref{eq:dL}) of the next section) and
$\delta\phi (\xi)=i_{\xi}\phi(\xi)$, where $\phi\in W^\s(\fg)$ is thought
of as a function on $\fg$ with values in $\wedge^\s\fg^*$. 
For example, when $\lambda\in\wedge^1\fg^*=\fg^*$, we have
$d_W\lambda=d_L\lambda-f_\lambda$ where $f_\lambda=\lambda$ is
viewed as an element of $S^1\fg^*$. The
differential $d_W$, the contractions $i_{\xi}$, and the natural action
of $G$ fit together to make $W^\s(\fg)$ into a $G$-differential algebra.
(See, e.g., \cite{fu:cohomology} or \cite{at-bo:moment_map} for more
details.) 

The Weil complex can be characterized by the following universal 
property.
(See, e.g., \cite{fu:cohomology} for the proof.) Namely, given a
multiplicative complex $A^\s$, any linear mapping 
$\Phi\colon\fg^*\to A^1$ can be extended in a unique 
way to a multiplicative homomorphism of complexes:
$\bar{\Phi}\colon W^\s(\fg)\to A^\s$. Note that 
$\bar{\Phi}(f_\lambda)=\wedge^2\bar{\Phi}(d_L\lambda)-d_A\Phi(\lambda)$,
where $\lambda$ and $f_\lambda$ are as above. The Weil complex is known 
to be acyclic. (See, e.g., \cite{fu:cohomology}, for an explicit homotopy
between $id$ and zero.)

Recall that for a $\fg$-differential complex $C^\s$, the basic   
subcomplex $C^\s_b$ is formed by $\fg$-invariant elements $c$ such    
that $i_{\xi}c=0$ for all $\xi\in\fg$. Equivalently, $c\in C^\s_b$    
if and only if $i_{\xi}c=0$ and $i_{\xi}dc=0$ for all $\xi$. The   
cohomology $H^\s_b(C^\s)=H^\s(C^\s_b)$ is called the basic cohomology   
of $C^\s$.

The tensor product of two $\fg$-differential complexes (algebras) is,  
in a natural way, a $\fg$-differential complex (algebra). In particular,   
given a $\fg$-differential complex (algebra) $A^\s$, the product   
$C^\s=A^\s\otimes W^\s(\fg)$ is a $\fg$-differential complex (algebra)   
and the basic cohomology $H^\s_b(A^\s\otimes W^\s(\fg))$ is defined.   

The inclusion $(S^\s\fg^*)^G\to (W^\s(\fg))_b$ is actually an isomorphism
because $i_\xi(\varphi\otimes f)=0$ for all $\xi$, where $\varphi\in\wedge^\s\fg^*$ and 
$f\in S^\s\fg^*$, means that $\varphi=1$. Observe also that
$d_W=0$ on $(S^\s\fg^*)^G$. Therefore,
$H^\s_b(W^\s(\fg))=(W^\s(\fg))_b=(S^\s\fg^*)^G$. In particular, this
implies that $H^\s_b(A^\s\otimes W^\s(\fg))$ is an $(S^\s\fg^*)^G$-module.

\begin{thm}
\label{thm:weil} 
$H^\s_b(A^\s\otimes W^\s(\fg))\simeq H^\s_G(A^\s)$ as topological
$(S^\s\fg^*)^G$-modules (or algebras if $A^\s$ is an algebra).
\end{thm}

The proof of the theorem can be found, for example, in 
Appendix B of \cite{gls}. The argument is rather standard: one shows
that the very complexes 
$(A^\s\otimes W^\s(\fg))_b$ and $(A^\s\otimes S^\s\fg^*)^G$ 
are isomorphic.\footnote{Given two
complexes, the result (e.g., Theorem \ref{thm:weil}) that their cohomology 
spaces are isomorphic may arguably be called ``trivial'' when the complexes 
are themselves isomorphic and ``nontrivial'' otherwise. This
classification does not appear to be completely meaningless.}
The isomorphism is just the natural inclusion
$(A^\s\otimes S^\s\fg^*)^G\to (A^\s\otimes W^\s(\fg))_b$ and its
inverse can be written explicitly.

\begin{dfn}
\label{dfn:locally-free}
{\em
A $\fg$-differential algebra $A^\s$ is said to be {\em locally free}
if $A^0$ has a unit 1, and so $\reals\subset A^1$, and there exists 
a $\fg$-equivariant linear homomorphism $\Theta\colon \fg^*\to A^1$ 
such that $\lambda(\xi)=i_\xi\Theta(\lambda)$ for all $\lambda\in \fg^*$
and $\xi\in\fg$.
}
\end{dfn}

\begin{rmk}
\label{rmk:average-free}
{\em 
If $G$ is compact, any $\Theta\colon \fg^*\to A^1$ 
satisfying the above condition $\lambda(\xi)=i_\xi\Theta(\lambda)$ can
be made equivariant by applying the averaging over $G$.

Definition \ref{dfn:locally-free} dates back to H. Cartan's 
original papers \cite{cartan1} and \cite{cartan2}. (Among more recent
publications see \cite{DKV} (Definition 16, the existence of
a connection) and also \cite{gls} (Appendix B).)
}
\end{rmk}

Consider now a $\fg$-differential algebra $A^\s$ and a   
$\fg$-differential complex $C^\s$. We say that $C^\s$ is   
a {\em $\fg$-differential $A^\s$-module} if $C^\s$ is 
a graded $A^\s$-module and the multiplication is compatible  
with the differentials and contractions. This means that 
$d(ac)=(da)c+(-1)^na\,dc$ and $i_\xi(ac)=(i_\xi a)c+(-1)^na(i_\xi c)$
for $a\in A^n$ and $c\in C^\s$.  

\begin{thm}  
\label{thm:basic-complex}  
Let $C^\s$ be a module over a locally free $\fg$-differential algebra.  
Then $H^\s_G(C^\s)\simeq H^\s_b(C^\s)$ as topological vector spaces.
\end{thm}

\begin{cor}
\label{thm:basic-algebra}  
Assume that $A^\s$ is locally free. Then $H^\s_G(A^\s)=H^\s_b(A^\s)$  as 
topological algebras.
\end{cor}

The proof of Theorem \ref{thm:basic-complex} is only a minor modification
of that of Corollary \ref{thm:basic-algebra}, which is essentially due
to H. Cartan. (See \cite{cartan1} and \cite{cartan2}.) The proofs can also
be found in \cite{DKV} (Theorem 17 of Part I, the ring of polynomials 
$S^\s\fg^*$ being replaced by  functions on $\fg^*$) and \cite{gls}
(Appendix B). We omit the argument here.

\begin{exam}
\label{exam:de-rham-loc-free}
{\em
(Continuing Example \ref{exam:de-rham}.) Let 
$G$ act on $M$ locally free, i.e., so that all stabilizers are
discrete. Then $A^\s=\Omega^\s(M)$ is a locally free $G$-differential
algebra. More generally, the de Rham complex is a locally free 
$\fg$-differential algebra when $M$ is given an infinitesimal
$G$-action generated by a monomorphism of $\fg$ into the space of
vector fields $\CX^1(M)$ on $M$. We shall call such an infinitesimal
action locally free as well. In the notation of Example \ref{exam:de-rham},
$\Omega^\s (M,E)$ is a module over $\Omega^\s (M)$, and the equivariant
cohomology is then the same as the cohomology of the basic complex
$\Omega^\s (M,E)_b$.

Assume now that $G$ is compact and that $M$ carries a genuine $G$-action.
Then, as well known, $H^\s_G(M)=H^\s_b(\Omega(M))=H^\s(M/G)$. 
One way to see this is to consider the Leray spectral sequence 
of the projection $\rho\colon (M\times EG)/G\to M/G$.
The fiber $\rho^{-1}(x)$, where $x\in M/G$, is homotopy
equivalent to $B G_x$. The Leray spectral sequence of $\rho$
with real coefficients collapses in the 
$E_2$-term, since $H^{\s>0}(B G_x,\reals)=0$. This yields an isomorphism
$H^\s_G(M)\simeq H^\s(M/G,\reals)$. (See, e.g., Lemma 5.3 of \cite{gi} 
and references therein
for more details.) Furthermore, the basic subcomplex $\Omega(M)^\s_b$ can be 
identified with the de Rham complex of the orbifold $M/G$. A similar
reasoning applies to the differential forms with values in a vector
bundle.

The Weil complex $W^\s(\fg)$ is only a nominally
different example of a locally free algebra.
Recall that geometrically $W^\s(\fg)$ can be thought of as an
algebraic model for $\Omega^\s(EG)$ and the $G$-action on $EG$ is free.
Clearly, $(W^\s(\fg))_b=(S^\s\fg^*)^G$ with the trivial differential,
and so $H^\s_b(W^\s(\fg))=(S^\s\fg^*)^G$.
}
\end{exam}

\subsection{The Chevalley--Eilenberg complex as a $K$-differential
complex}
\label{subsec:chev-eilen}

\subsubsection{The Chevalley--Eilenberg complex}

Let $K\subset G$ be a connected subgroup of $G$ with the Lie algebra $\fk$
and let $V$ be a locally convex differentiable
$G$-module. Consider the complex 
$C^\s(\fg,K;V)=(V\otimes\wedge^\s(\fg/\fk)^*)^K$ with the differential
$d_{Lie}$ defined as
\begin{eqnarray}   
\label{eq:dL}   
d_{Lie}\phi(\xi_1,\ldots,\xi_{n+1})&=&   
\sum_{1\leq j<j\leq n+1}    
(-1)^{j+j-1}\phi([\xi_i,\xi_j],   
\xi_1,\ldots,\hat{\xi}_i,\ldots,\hat{\xi}_j,\ldots,\xi_{n+1})\\    
\nonumber &&\quad+\sum_{1\leq k\leq n+1}(-1)^{k}    
\xi_k\phi(\xi_1,\ldots,\hat{\xi}_k,\ldots,\xi_{n+1})  
\enspace ,    
\end{eqnarray}
where $\xi v$ denotes the action of $\xi\in \fg$ on $v\in V$ and
$C^\s(\fg,K;V)$ is identified with a subspace of 
$C^\s(\fg;V)=V\otimes\wedge^\s\fg^*$. In this way, we do obtain
a complex  called the 
Chevalley--Eilenberg complex of $\fg$ relative to $K$ with coefficients
in $V$. (See, e.g., \cite{fu:cohomology} or\cite{gui:book} for more details.)
The cohomology of this complex is denoted by $H^\s(\fg,K;V)$.

Take, in particular, $(C^\s,d)=(C^\s(\fg,V),d_{Lie})$. A straightforward
calculation shows that the natural $G$-action 
and contractions make $C^\s$ into a $G$-differential complex,
and the equivariant cohomology $H^\s_K(C^\s)$ is defined.
Note also that if $V$ is a topological $G$-algebra, i.e., the multiplication
$G$-invariant, $C^\s(\fg,K;V)$ is a multiplicative complex, and so
$C^\s$ is a $G$-differential algebra. In any case

\begin{thm}
\label{thm:eq-Lie}
Assume that $\fk\subset\fg$ has a $K$-invariant complement with respect
to the adjoint action. Then $H^\s_K(C^\s)=H^\s(\fg,K;V)$ as topological
vector spaces. Furthermore, 
these spaces are isomorphic as topological algebras if $V$ is a $G$-algebra.
\end{thm}

\begin{rmk} 
{\em 
A $K$-invariant complement $\fk^{\perp}$ exists, for example, when
$K$ is compact or $\fk$ is semisimple. (But also in some other cases as well.)
The existence of $\fk^{\perp}$ is equivalent to the existence of
a $K$-equivariant projection $\theta\colon\fg\to\fk$, for we can take
$\fk^{\perp}=\ker\theta$.
}
\end{rmk}

\proof Observe that $A^\s=C^\s(\fg)$ is a free $\fk$-differential algebra. 
(The linear map $\Theta=\theta^*\colon \fk^*\to \fg^*$ is $\fk$-equivariant.)
Furthermore, $C^\s=C^\s(\fg;V)$ is a differential $A^\s$-module. Hence
by Theorem \ref{thm:basic-complex}, $H^\s_K(C^\s)=H^\s_b(C^\s)$.
It is easy to see that $C^\s(\fg;V)_b=C^\s(\fg,K;V)$ as complexes,
and $H^\s_b(C^\s)=H^\s(\fg,K;V)$.\footnote{According to the
classification of our previous footnote, Theorem \ref{thm:eq-Lie},
as well as Theorem \ref{thm:basic-complex}
it depends on, is ``nontrivial''.}
\qed

\begin{rmk}
{\em
It is elucidating to discuss the geometrical meaning of Theorem
\ref{thm:eq-Lie} in the context of Example \ref{exam:de-rham}.
Consider the $K$-principal bundle $K\to G\stackrel{\pi}{\to} G/K$. 
Let $E=G\times V$ and $E'=(G/K)\times V$ be the trivial bundles equipped with 
the canonical flat connections. As in Example \ref{exam:de-rham},
$\Omega^\s(G,V)$ is a differential $K$-complex with respect to the right
action of $K$ on $G$. Since the action is free, we have
$H^{\s}_K(\Omega^\s(G;V))=H^\s(G/K;V)$. The isomorphism is
induced by the pull-back
$\pi^*\colon \Omega^\s(G/K;V)\to \Omega^\s(G;V)_b$ (cf., \cite{cartan2}).
The complex $C^\s(\fg,K;V)$ can be identified with the subcomplex
of $\Omega^\s(G/K;V)$ formed by $V$-valued differential forms invariant
under the diagonal action. The averaging yields 
$H^\s(\fg,K;V)=H^\s_{dR}(G/K;V)$, provided that $G$ is compact. (See
the next section for a detailed discussion.)  Similarly,
$C^\s(\fg;V)$ is a subcomplex of forms invariant under
the diagonal action on $\Omega^\s(G;V)$. It is not hard to see that
$C^\s(\fg;V)$ and $\Omega^\s(G;V)$ have equal $K$-equivariant 
cohomology when $G$ is compact. Thus
$$
H^\s(\fg,K;V)=H^\s(G/K;V)=H^{\s}_K(\Omega^\s(G;V))=H^\s_K(C^\s)
\enspace ,$$
which gives a topological proof of Theorem \ref{thm:eq-Lie} for
a compact $G$.
}
\end{rmk}

\subsubsection{Lie algebra cohomology with coefficients in Fr\'{e}chet 
modules}
\label{section:lie-alg}

We have just seen how the explicit calculation of certain equivariant 
cohomology spaces can be reduced to finding relative Lie algebra 
cohomology. In the subsequent sections 
we shall see that some Poisson cohomology can also be treated as the 
Lie algebra cohomology. 
Various Lie algebra or Lie group cohomology spaces with coefficients in 
other Fr\'{e}chet modules also arise in studying the rigidity of group
actions (e.g., Poisson actions) and the existence of pre-momentum
mappings \cite{gi:moment}. In this section we recall how to
carry out the calculation of cohomology explicitly when the group is compact.

Let $G$ be a compact Lie group and $V$ a smooth Fr\'{e}chat $G$-module.
Denote by $\rho$ the representation on $G$ on $V$ and by $V^G$
the (closed) subspace of all $G$-invariant vectors in $V$.
Recall that for every $v\in V$, the map $g\mapsto \rho(g)v$ of $G$ 
to $V$ is smooth and therefore $V$ is a $\fg$-module. 
Furthermore, this map is integrable with respect to a bi-invariant Haar 
measure $dg$ on $G$ and the projector $A\colon V\to V^G$ given by averaging
over $G$ is continuous. 
Let $C^\s$ be the the Chevalley--Eilenberg complex
$C^{\s}(\fg, K; V)$ of $\fg$ relative $K$ with values in $V$. 
In this section we use the notation $d_{\rho}$ for the differential 
$d_{Lie}$ of (\ref{eq:dL})
to emphasize its dependence of $\rho$. Denote $B^n\subset C^n$ the 
subspace of all exact $n$-cochains and by $Z^n$ the subspace formed by
cocycles. The space $C^n$ is Fr\'{e}chet as well as $Z^n$, and $B^n$,
being equipped with the induced topology, turns into a locally convex 
topological vector space.

Sometimes one also needs to consider a family $\rho_x$, $x\in X$, 
of representations of $G$ on $V$ smoothly parameterized by a manifold $X$.
Then the operator $d_{\rho_x}$ depends on $x$ but the spaces $C^n$ do not.
 
\begin{thm}
\label{thm:Lie-cohomology} 
~~~
\begin{enumerate}

\item[{\rm (i)}] The inclusion $V^G\to V$ gives rise to a topological 
isomorphism in the
Lie algebra cohomology whose inverse is induced by $A$. In particular,
$$
H^{\s}(\fg, K; V)=H^{\s}(\fg, K; V^G)=H^{\s}(\fg, K)\otimes V^G
\enspace .$$

\item[{\rm (ii)}] For every $n>0$, there exists a continuous linear map 
$$
\CH\colon C^n(\fg, K; V)\longrightarrow C^{n-1}(\fg, K; V)
$$
such that $d\CH\mid_{B^n}=id$.

\item[{\rm (iii)}] Given a smooth family of representations $\rho_x$, 
the family of operators $\CH_x$ from {\rm (ii)} can be chosen to be 
smooth in $x$. 
\end{enumerate}
\end{thm}

\begin{cor}  
\label{cor:equi-expl}  
In the notation of Theorem {\em \ref{thm:eq-Lie}},   
$H^{\s}_K(C^\s)=H^\s(\fg,K)\otimes V^G$, provided that $G$ is compact and    
$V$ is a Fr\'{e}chet $G$-module.  
\end{cor}

Assertion (i), which is sufficient for a majority of applications,
can be easily proved
by employing the so-called  van Est spectral sequence \cite{es:application},
e.g., as in \cite{gi:moment}, combined with a theorem by Mostow which
implies that $\Omega^\s(G/K;V)$ is continuously strongly injective.
(See \cite{bw}, Proposition 5.4 of Chapter IX.)
In the appendix, we outline a direct and relatively
elementary proof, yet also using Mostow's theorem, of all the assertions.

Another immediate consequence of the theorem is the following 

\begin{cor}   
\label{cor:coh-functions}  
Let a compact group $G$ act on a manifold $M$. Then   
$H^{\s}(\fg; C^{\infty}(M))=H^{\s}(\fg)\otimes (C^{\infty}(M))^G$.
\end{cor}

\begin{rmk}
{\em
1. For a finite-dimensional $V$,
assertion (i) of Theorem \ref{thm:Lie-cohomology}
is well known (e.g., \cite{gui:book}, Section II.11) and (ii), (iii)
are evident. Moreover, it is a classical result that $H^{\s}(\fg; V)=0$, 
when $\fg$ is a semisimple real or complex finite-dimensional Lie
algebra and $V$ an irreducible finite-dimensional nontrivial
$\fg$-module. (See, e.g., \cite{gui:book} or \cite{fu:cohomology} and 
references therein.) Of course, this also follows from (i) by complexifying, 
in the real case, both $\fg$ and $V$ and then passing to the representation
of the compact form of $\fg_{\complex}$ on $V_{\complex}$. 

2. 
It is essential that the $\fg$-module structure on $V$ integrates to that
of a $G$-module. Theorem \ref{thm:Lie-cohomology} does not extend to
just $\fg$-modules \cite{gw}. 
Nor does it hold when $G$ is a noncompact
simple group (e.g., $G={\rm SL}(2,\reals)$).
}
\end{rmk}

\subsection{Further properties of equivariant cohomology}
\label{subsec:further}

In this section, we routinely extend the standard properties
of the equivariant de Rham cohomology to the equivariant cohomology of 
$G$-differential complexes.

Let us start with a calculation of equivariant cohomology
in low degrees. The zeroth cohomology is particularly simple:
$H^0_G(A^\s)=\ker(d\colon A^0\to A^1)$,
for the elements of the kernel are $\fg$-invariant.

\begin{thm} ~~~~
\label{thm:low}
\begin{enumerate} 

\item[{\em (i)}] Let $H^1(\fg)=0$. Then every $d$-closed $\fg$-invariant 
element of $A^1$ is automatically $d_{\fg}$-closed, and, as a consequence,
$$
H^1_{G}(A)
=\{a\in A^1\mid~d a=0~\hbox{\rm and}~ i_{\xi}a=0
~\hbox{\rm for all}~\xi\in\fg\}/\,d(A^0)^G
\enspace .$$ 

\item[{\em (ii)}] Let $G$ be compact semisimple and let $A^\s$ be
a $G$-differential complex. Then the forgetful
homomorphism $H^\s_{G}(A)\to H^\s(A)$ is an isomorphism in degree one,
$\s=1$, and an epimorphism in degree two, $\s=2$.
\end{enumerate}
\end{thm}

\proof Denote by $Z^n$ the space of all $d$-closed elements in $A^n$. It is
clear that $Z^n$ is a Fr\'{e}chet $\fg$- or $G$-submodule of $A^n$.
We will need the following observation.
\begin{itemize}
\item {\em Let $a\in (Z^n)^\fg$. Then $c\colon\fg\to Z^{n-1}$, defined as 
$c(\xi)=i_\xi a$, is a $\fg$-cocycle.}
\end{itemize}
This fact an immediate consequence of the identity
$i_{[\xi,\zeta]}a=L_\xi i_\zeta a-i_\zeta L_\xi a$ combined with
the assumption that $a$ is closed and $\fg$-invariant.

Let us prove (i). We have to show that 
$$
 i_{\xi}a=0~\hbox{\rm for all}~ \xi\in\fg \Longleftrightarrow
L_{\xi}a=0~\hbox{\rm for all}~ \xi\in\fg \enspace ,
$$
when $a\in Z^1$. Clearly, $da=0$ and $i_{\xi}a=0$ yield $L_{\xi}a=0$.
Conversely, assume that the right hand condition holds. Then $c$
defined as above is a cocycle on $\fg$ with coefficients in $Z^0$.
Since $Z^0$ is a trivial $\fg$-module, $H^1(\fg;Z^0)=H^1(\fg)\otimes Z^0=0$, 
and $c$ is exact. An exact one-cocycle with coefficients in a trivial 
module is identically zero. Hence $i_{\xi}a=0$ for all $\xi\in\fg$.

The forgetful homomorphism is injective in degree one because
$A^0_\fg=(A^0)^G$. Hence we just need to show that it is onto.
Let $a\in A^1$ be $d$-closed. Without loss of generality we
may assume that $a$ is $G$-invariant. (Indeed, since $G$ is
compact, the mean $\int_G ga\,dg$ exists, is $G$-invariant 
and homologous to $a$ due to Cartan's formula.) Then by (i),
$a$ is also $d_\fg$-closed.

It remains to show that the forgetful homomorphism is surjective in degree two.
Pick $a\in Z^2$. As before, we may assume that $a$ is $G$-invariant, i.e.,
$a\in (Z^2)^G$. Then $c(\xi)=i_\xi a$ is a cocycle on $\fg$ with values
in $Z^1$. Since $\fg$ is semisimple, $H^1(\fg)=0$. Thus, by Theorem 
\ref{thm:Lie-cohomology}, $H^1(\fg; Z^1)=H^1(\fg)\otimes Z^1=0$, and
$c$ is exact. In other words, there exists 
$b\in (Z^1)^G$ such that $i_\xi a=L_\xi b$. It is clear that the element 
$a+i_{\bullet}b\in A^2_\fg$ is $d_\fg$-closed and projects onto $a$ under
the forgetful homomorphism. \qed

Let us now turn to some general properties of equivariant cohomology.
The following proposition is evident:  
\begin{prop}   
{\em (The long exact sequence.)}   
Let $0\to B^\s\to A^\s\to C^\s\to 0$ be an exact sequence of   
$G$-differential complexes. If $G$ is compact, the induced sequence   
$$
0\longrightarrow B^\s_\fg\longrightarrow A^\s_\fg\ \longrightarrow C^\s_\fg\longrightarrow 0 
$$ 
is also exact, and we have a long exact sequence for equivariant cohomology:   
$$   
\ldots\longrightarrow H^\s_G(B^\s)\longrightarrow H^\s_G(A^\s)\longrightarrow H^\s_G(C^\s)\longrightarrow \ldots   
\enspace .$$   
\end{prop}

\begin{rmk}
{\em
Let $A^\s=\Omega^\s(M)$, where $M$ is a smooth manifold. The
splitting map $p\colon (M\times E\T)/\T\to (M\times EG)/G$ is homotopy
equivalent to a fibration with fiber $G/\T$. The map $p$ is known to 
induce a monomorphism $H^{\s}_G(M)\to H^\s_\T(M)$ whose image is
$H^\s_\T(M)^W$, where $W$ is the Weil group of $G$. This result can also 
be generalized to equivariant cohomology of $G$-differential complexes when 
$G$ is compact (cf., \cite{DKV}, p. 155). We omit the details since this 
fact, although very interesting, is irrelevant to the subject matter of 
this paper.
} 
\end{rmk}

To introduce the standard spectral sequence for equivariant cohomology,
consider the decreasing filtration 
$$
F_0=A^{\s}_G\supset F_1\supset\ldots\supset F_n\supset F_{n+1}\supset\ldots
$$
of the complex $A^{\s}_G$ by its subcomplexes
$$
F_n=\bigoplus_{2j\geq n}(A^\s\otimes S^j\fg^*)^G
\enspace . $$

\begin{thm}
\label{thm:spec-seq}
The spectral sequence of the filtration $F_p$ converges to
$H^{\s}_G(A^\s)$. If $G$ is compact,
$$
E^{pq}_1=E_2^{pq}=H^{q}(A^\s)\otimes (S^{p/2}\fg^*)^G
$$
when $p$ is even and $E^{pq}_1=E_2^{pq}=0$ when $p$ is odd. 
Furthermore,
$$
d_2([a]\otimes\phi)(\xi)=[i_{\xi}a]\phi(\xi)
\enspace, $$
where $\phi\in S^\s\fg^*$ and $[a]$ is the cohomology
class of $a\in A^\s$.
\end{thm}

\begin{rmk}
{\em
The $E_2$-term can be sometimes found
even when $G$ is not compact. For example, 
if $G$ is abelian, we have $E_2^{pq}=H^{q}((A^\s)^G)\otimes S^{p/2}\fg^*$.
(Generically zero!)
It is clear that in general $E_2\neq H^{\s}(A^\s)\otimes (S^\s\fg^*)^G$. 

Theorem \ref{thm:low} can be derived from Theorem \ref{thm:spec-seq} .
However, a direct elementary proof given above seems more transparent.

When $A^\s=\Omega^\s(M)$, the spectral sequence turns into the Leray
spectral sequence of the fibration 
$M\to (M\times EG)/G\to BG$.\footnote{As a matter of tradition, when the
$E_2$-term of the Leray spectral sequence of a fibration is expressed
as the tensor product, the cohomology of the base is taken first and
the cohomology of the fiber second. This convention is systematically 
broken throughout the paper. To follow the rule, one would have to 
switch the order in the tensor product, i.e., to start with 
$(S^\s\fg^*\otimes A^\s)^G$, so that $E_2$ would become
$(S^\s\fg^*)^G\otimes H^\s(A^\s)$. Of course, the results obtained for
different orders are canonically isomorphic. The same applies to the 
spectral sequence introduced in the next section.}
}
\end{rmk}

\proof The convergence  of $E_r$ is a general fact which
holds for any spectral sequence arising from a filtered complex.

To evaluate the $E_2$-term, note that by definition
$E_1^{pq}$ is the $(p+q)$-th cohomology of $E_0^{p\s}=F_p^\s/F_{p+1}^\s$ 
with respect to the induced differential $d_0$. It is clear that $E_0^{p\s}$ 
is just $(A^{\s}\otimes S^{p/2}\fg^*)^G$ when $p$ is even and zero 
otherwise. The differential  $d_0$ is the restriction of the original
differential $d$ of $A^\s$ extended to act trivially on $S^{p/2}\fg^*$.
When $G$ is compact, we have 
$$
H^\s((A^{\s}\otimes S^{p/2}\fg^*)^G)=
H^\s(A^{\s}\otimes S^{p/2}\fg^*)^G
\enspace .$$
By the definition of $d_1$, the right hand side is just
$(H^\s(A^{\s})\otimes S^{p/2}\fg^*)^G$.
Furthermore, since $A^\s$ is a differential $G$-module, $G$ acts trivially
on the cohomology of $A^\s$. Thus
$E^{pq}_1=H^q(A^{\s})\otimes (S^{p/2}\fg^*)^G$. Clearly, $d_1=0$.
Hence, $E_2=E_1$ and $d_2$ is simply induced by the differential $d_G$.
\qed

\subsection{The second spectral sequence for equivariant cohomology}
\label{sec:spec-general}

In this section we introduce a spectral sequence relating the equivariant
cohomology of a differential $G$-complex $A^\s$, the cohomology of $G$,
and the cohomology of $A^\s$. To be more precise, the spectral sequence 
converges to $H^\s(A^\s)$ and has the product 
$H^\s(\fg)\otimes H^\s_G(A^\s)$ 
as its $E_2$-term. This result will be used in Section \ref{sec:spec-mom} 
to construct a spectral sequence in Poisson cohomology associated
with an equivariant momentum mapping.

\begin{exam}
\label{exam:spec-man}
{\em
Let $M$ be a smooth manifold acted on by a compact group $G$. 
As in Example \ref{exam:de-rham}, take $A^\s=\Omega^\s(M)$. If the action
is free, the natural projection $M\times EG\to (M\times EG)/G$ 
is a principal $G$-bundle. The Leray spectral sequence of this bundle 
converges to $H^\s(M)$ because $EG$ is contractible and has 
$E_2=H^\s(G)\otimes H^\s_G(M)$. The same is true when the action is just
locally free.
}
\end{exam}

\begin{thm}
\label{thm:spec-new}
Assume that $G$ is compact.
There exists a spectral sequence converging
to $H^\s(A)$ and such that 
$E_1^{pq}=H^q(\fg)\otimes (A^\s\otimes W^\s(\fg))^p_b$ and
$E^{pq}_2=H^q(\fg)\otimes H^p_G(A^\s)$.
\end{thm}

\begin{rmk}
{\em
As clear from the proof, when $G$ is not compact, one still has
$E_1^{pq}=H^q(\fg; (A^\s\otimes W^\s(\fg))^p_0)$, where
$(A^\s\otimes W^\s(\fg))_0$ is the submodule of $A^\s\otimes W^\s(\fg)$
formed by $a\otimes w$ which have zero contraction, $i_\xi(a\otimes w)=0$,
with every $\xi\in\fg$.
}
\end{rmk}

\proof The $G$-differential complex $B^\s=A^\s\otimes W^\s(\fg)$ 
carries the following decreasing filtration by subcomplexes:
$$
F^{p+q}_p=\{ b\in B^{p+q}\mid i_{\xi_1}i_{\xi_2}\ldots i_{\xi_{q+1}}b=0
~\hbox{\rm for all}~\xi_1,\ldots, \xi_{q+1}\in \fg\}
\enspace .$$
In degree $n$, the smallest nonzero term $F^n_n$ of this filtration is
formed by $b\in B^n$ such that $i_\xi b=0$ for all $\xi\in\fg$.
The spectral sequence $E_r$ in question is just that associated with the
filtration. It is clear that $E_r$ converges to $H^\s(A)$ because $W^\s(\fg)$
is acyclic.

We claim that $E^{pq}_0=F^{p+q}_p/F^{p+q}_{p+1}$ is naturally isomorphic to 
$C^q(\fg; F^p_p)=\hbox{\rm Hom}(\wedge^q\fg, F^p_p)$. 
For $b\in F^{p+q}_p$, let 
$$
c_b(\xi_1,\ldots,\xi_q)=i_{\xi_q}\ldots i_{\xi_1}b
\enspace. $$
By the definition of the filtration, $c_b(\xi_1,\ldots,\xi_q)\in F^p_p$.
Thus $c_b$ can be viewed as a liner mapping $\wedge^q\fg\to F^p_p$. 
Furthermore, $c$ is also linear in $b$ and $c_b\equiv 0$ if and only if 
$b\in F^{p+q}_{p+1}$. This means that $b\mapsto c_b$ is a monomorphism
$E^{pq}_0\to C^q(\fg; F^p_p)$. This
monomorphism is actually an epimorphism, and so an isomorphism. To
see this, let us pick $\phi\in C^q(\fg; F^p_p)$ 
and show that $\phi=c_b$ for some $b\in F^{p+q}_p$. 
In effect, we may just take $b=\sum \phi(e_j)\otimes e^*_j$, where
$e_1,\ldots, e_r$ is a basis in $\wedge^q\fg$ and
$e_1^*,\ldots, e_r^*$ is the dual basis.

The next step is to prove that the differential $d_0$ is, up to a sign,
the Chevalley--Eilenberg differential $d_{Lie}$ on $C^q(\fg;F^p_p)$ with
$F^p_p$ equipped with its natural $\fg$-module structure inherited 
from $B^\s$. (Note in this connection that $F^p_p$ is indeed a submodule 
of $B^p$ as follows from (ii'): 
$i_{[\xi,\zeta]}=L_\xi i_\zeta- i_\zeta L_\xi$.) This is equivalent to showing
that $c_{db}=-d_{Lie}c_b$, which is easy to verify by a straightforward 
calculation. Namely, using Cartan's formula, we obtain
\begin{eqnarray*}
c_{db}(\xi_1,\ldots,\xi_{q+1})&=&i_{\xi_{q+1}}\ldots i_{\xi_1}db\\
&=&-i_{\xi_{q+1}}\ldots i_{\xi_2} d i_{\xi_1}b+i_{\xi_{q+1}}\ldots L_{\xi_1}db\\
&=&i_{\xi_{q+1}}\ldots i_{\xi_3} di_{\xi_2}i_{\xi_1}b 
-i_{\xi_{q+1}}\ldots i_{\xi_3} L_{\xi_2}i_{\xi_1}b 
+ i_{\xi_{q+1}}\ldots L_{\xi_1}db\\
&\vdots&\quad \\
&=&\sum_{1\leq k\leq q+1}(-1)^{k-1}
i_{\xi_{q+1}}\ldots i_{\xi_{k+1}} L_{\xi_k}i_{\xi_{k-1}}\ldots i_{\xi_1}b 
+(-1)^{q+1}d i_{\xi_{q+1}}\ldots i_{\xi_1}b\\
&=&\sum_{1\leq k\leq q+1}(-1)^{k-1}
i_{\xi_{q+1}}\ldots i_{\xi_{k+1}} L_{\xi_k}i_{\xi_{k-1}}\ldots i_{\xi_1}b 
\enspace ,
\end{eqnarray*}
where the last term vanishes because $b\in F^{p+q}_p$. Employing (ii'), we
can switch $L_{\xi_k}$ and the contractions so that
to have $L_{\xi_k}$ applied last:
\begin{eqnarray*}
c_{db}(\xi_1,\ldots,\xi_{q+1})&=&
-\sum_{1\leq j<l \leq q+1}
(-1)^{j-1}
i_{\xi_{q+1}}\ldots i_{\xi_{l+1}} i_{[\xi_j,\xi_l]} i_{\xi_{l-1}}
\ldots \widehat{i_{\xi_j}}\ldots i_{\xi_1} b\\
&&\quad+\sum_{1\leq k \leq q+1} (-1)^{k-1}
L_{\xi_k}i_{\xi_{q+1}}\ldots \widehat{i_{\xi_k}}\ldots i_{\xi_1} b\\
&=&
\sum_{1\leq j<l \leq q+1}(-1)^{l+j}
i_{\xi_{q+1}}\ldots \widehat{i_{\xi_{l}}}\ldots \widehat{i_{\xi_{j}}}\ldots 
i_{\xi_{1}}\ldots i_{[\xi_{j},\xi{l}]} b \\
&&\quad + \sum_{1\leq k \leq q+1} (-1)^{k-1}
L_{\xi_k}i_{\xi_{q+1}}\ldots \widehat{i_{\xi_{k}}}\ldots i_{\xi_{1}} b\\
&=&-d_{Lie}c_b(\xi_1,\ldots,\xi_{q+1})
\enspace .
\end{eqnarray*}

The $E_1$ term is now easy to find: $E^{pq}_1=H^q(\fg; F^p_p)$. By
Theorem \ref{thm:Lie-cohomology}, this is just $H^q(\fg)\otimes (F^p_p)^G$
because $G$ is compact and $F^p_p$ is clearly a Fr\'{e}chet $G$-module. 
Finally, recall that $F^p_p$ is formed by $b\in B^p$ such that 
$i_\xi b=0$ for any $\xi\in\fg$. Thus $(F^p_p)^G= B^p_b$, which yields
$E_1^{pq}=H^q(\fg)\otimes (A^\s\otimes W^\s(\fg))^p_b$.

It remains to prove that $d_1$ is just the natural differential $d$ 
inherited by $B^\s_b$ from $B^\s$ and acting on the second term 
only in the tensor product expression for $E^{pq}_1$. To this end, pick a 
representative $\phi\otimes f$ of an element of $E^{pq}_1$ such that
$\phi$ is a $q$-cocycle on $\fg$ and $f\in B^p_b$. We need to show
that $d_1(\phi\otimes f)=\phi\otimes df$. Recall that 
$d_1$ is induced by $d= d_A\pm d_W$, where $d_W=d_L-\delta$. (See Section
\ref{sec:Weil}.) Since $\phi$ is a cocycle, $d_L\phi=0$, and so
$d(\phi\otimes f)=\pm (\delta\phi)\otimes f+\phi\otimes df$. The first
term $(\delta\phi)\otimes f$ belongs to $F^{p+q+1}_{p+2}$, because $f$
is basic. Thus it makes no contribution into $d_1$, which therefore
has the desired form. As a result,
$E^{pq}_2=H^q(\fg)\otimes H^p_b(B^\s)=H^q(\fg)\otimes H^p_G(A^\s)$.
\qed

One can also start with a filtration of $A^\s$ instead of the filtration
of $A^\s\otimes W^\s$ considered above. Namely, let us take
\begin{equation}
\label{eq:filt}
F^{p+q}_p=\{ a\in A^{p+q}\mid i_{\xi_1}i_{\xi_2}\ldots i_{\xi_{q+1}}a=0
~\hbox{\rm for all}~\xi_1,\ldots, \xi_{q+1}\in \fg\}
\enspace .
\end{equation}
This filtration gives
rise to a spectral sequence converging to $H^\s(A^\s)$ but
having in general an absolutely incomprehensible structure.
However, if $A^\s$ is locally free,  the spectral sequence is
even simpler than that of Theorem \ref{thm:spec-new} and
essentially equivalent to it.

\begin{exam}
\label{exam:spec-man-free}
{\em
In the notation of Example  \ref{exam:spec-man}, the projection 
$M\to M/G$ gives rise to a filtration on $\Omega^\s(M)$ but rather little
can be said about the resulting spectral sequence  except that it converges 
to $H^\s(M)$ when no extra assumption on the action is made.
However, if the action is locally
free, we have $E^{pq}_1=H^q(G)\otimes \Omega^p_b(M)$ and
$E_2^{pq}=H^{q}(G)\times H^p(M/G)$. In fact, the spectral sequence
of this projection coincides with that from Example \ref{exam:spec-man}
beginning with the $E_2$-term. When the action is free, 
the fibration $M\to M/G$ is homotopy equivalent to 
$M\times EG\to (M\times EG)/G$.
}
\end{exam}

The proof of our next theorem is omitted, for it is quite similar
to that of Theorem \ref{thm:spec-new} (cf., the proof of Theorem 
\ref{thm:spec-poisson-sub}).

\begin{thm}
\label{thm:spec-new-free}
Assume that $G$ is compact and $A^\s$ is a locally free $G$-differential
algebra. The spectral sequence associated with the filtration 
{\em (\ref{eq:filt})} converges to $H^\s(A^\s)$ and has 
$E^{pq}_1=H^q(\fg)\otimes A^p_b$ and $E^{pq}_2=H^q(\fg)\otimes H^p_b(A^\s)$.
\end{thm}

\begin{rmk}
{\em
In fact, Theorem \ref{thm:spec-new} follows
from Theorem \ref{thm:spec-new-free} under a very natural extra 
hypothesis that $A^\s$ is a $G$-differential algebra with unit. Indeed,
then $B^\s=A^\s\otimes W^\s(\fg)$ becomes a locally free $G$-differential
algebra, and applying Theorem \ref{thm:spec-new-free} to $B^\s$, we obtain 
Theorem \ref{thm:spec-new}. Furthermore, it is not hard to extend Theorem
\ref{thm:spec-new-free} to complexes over locally free $G$-differential
algebras. Then it would yield Theorem \ref{thm:spec-new} in its full 
generality.
}
\end{rmk}

\section{Some results from Poisson geometry} 
\label{sec:dfn} 

This section is a brief review on the geometry of Poisson manifolds
heavily biased toward our present needs. 
In addition to the standard material, which can be found, for example, in
\cite{va:book}, we also prove some new theorems and recall the results 
of \cite{gi:moment} essential for what follows.

\subsection{Poisson manifolds and momentum mappings}

Let $P$ be a Poisson manifold. Denote by $\CXP$ the algebra of multi-vector
fields on $P$.
Recall that the space of one-forms $\Omega^1(P)$ on $P$ is a (local) Lie algebra with the bracket   
\begin{eqnarray}   
\label{eq:bracket}   
\{\alpha,\beta\}&=&d\left<\alpha\wedge\beta,\pi\right>    
+i_{\pi^{\#}\alpha}d\beta-i_{\pi^{\#}\beta}d\alpha\\    
\nonumber &=&+L_{\pi^{\#}\alpha}\beta-i_{\pi^{\#}\beta}d\alpha   
\enspace ,    
\end{eqnarray}   
where $\alpha$ and $\beta$ are one-forms, $\left<~,~\right>$ denotes   
the pairing between $k$-forms and $k$-vector fields and   
$\pi^{\#}\colon \Omega^1(P)\to\CX^1(P)$ is the paring with $\pi$, i.e.,
$\beta(\pi^{\#}\alpha)=\left<\beta\wedge\alpha,\pi\right>$.   
Note that $\pi^{\#}$ and $d\colon C^{\infty}(P)\to \Omega^1(P)$ are 
Lie algebra homomorphisms and $\pi^\#$ extends to an algebra homomorphism
$\Omega^\s(P)\to\CX^\s(P)$ by multiplicativity: 
$\pi^\#(\alpha_1\wedge\alpha_2)=\pi^\#(\alpha_1)\wedge\pi^\#(\alpha_2)$.

Following Bhaskara and Viswanath \cite{bv}, let us define
the ``Lie derivative'' of a multi-vector field $w$ in the direction
of a one-form $\alpha$ by Cartan's identity:
\begin{equation}   
\label{eq:Cartan-Poisson}    
\CL_{\alpha}w=i_{\alpha}d_{\pi}w + d_{\pi}i_{\alpha}w   
\enspace, \end{equation}  
where $d_\pi w=-[\pi,w]$, with $[~,~]$ being the Schouten bracket, and
$i_{\alpha}w$ stands for the contraction of $\alpha$ with   
$w$. To be more precise, $i_{\alpha}w$ is characterized by the equation 
$\left<\beta, i_{\alpha}w\right>=\left<\alpha\wedge\beta,w\right>$.  

As shown in \cite{bv}, the action $\CL$ makes the space of multi-vector 
fields $\CX^\s(P)$ into a module over $\DIF(P)$:
$$
\CL_{\{\alpha,\beta\}}=\CL_{\alpha}\CL_{\beta}-\CL_{\beta}\CL_{\alpha}
$$   
and 
\begin{equation}   
\label{eq:ii}   
i_{\{\alpha,\beta\}}=\CL_{\alpha}i_{\beta}-i_{\beta}\CL_{\alpha}   
\enspace .   
\end{equation}   
Furthermore, $\CL_{\alpha}$ satisfies the Leibniz identity
\begin{equation}   
\label{eq:Leibn}   
\CL_{\alpha}(w_1\wedge w_2)=(\CL_{\alpha}w_1)\wedge w_2+   
w_1\wedge(\CL_{\alpha}w_2)  
\enspace ,\end{equation}  
where $w_1$ and $w_2$ are multi-vector fields. In accord with
(\ref{eq:Cartan-Poisson}) we have
$$   
\CL_{\alpha}f=L_{\pi^{\#}\alpha}f  
$$   
for $f\in C^{\infty}(P)$. The ``Lie derivative'' $\CL_{\alpha}$ 
of a vector field $v$ is
related to the standard Lie derivative along $\xi=\pi^{\#}\alpha$ by 
the following equation verified in \cite{gi:moment}:
\begin{equation}   
\label{eq:Poisson-Lie-1}   
\CL_\alpha v   =L_\xi v+\pi^\# i_{v}\, d\alpha
\enspace  .
\end{equation}
As a consequence,
\begin{equation}   
\label{eq:Poisson-Lie-2}   
\CL_\alpha w   =L_\xi w+
\sum_{j=1}^q   
v_1\wedge\ldots\wedge v_{j-1}\wedge(\pi^{\#}i_{v_j}d\alpha)   
\wedge v_{j+1}\wedge \ldots\wedge v_q   
\enspace ,  \end{equation}
where $v_1,v_2,\ldots,v_q$ are vector fields and 
$w=v_1\wedge v_2\wedge\ldots\wedge v_q$ \cite{gi:moment}. This shows that
$\pi^{\#}\colon \OP\to\CXP$ is a homomorphism of $\DIF(P)$-modules,
provided that the domain is given the $\DIF(P)$-action via the standard
Lie derivative $L$ while the range via $\CL$. In other words,  
\begin{equation}   
\label{eq:alt-mod}   
\pi^{\#}L_{\pi^{\#}\alpha}\beta=\CL_{\alpha}\pi^{\#}\beta   
\enspace ,  
\end{equation}  where $\alpha\in\DIF(P)$ and $\beta\in\OP$. 

With the forthcoming applications in mind, let us
rewrite (\ref{eq:Poisson-Lie-2}) in a more compact form. 
For $w\in \CX^q(P)$ consider a $C^{\infty}(P)$-linear homomorphism
$$
\tilde{i}_w\colon \Omega^p(P)\longrightarrow\Omega^{p-1}(P)\otimes \CX^{q-1}(P)
$$
defined by the formula
\begin{equation}   
\label{eq:ti}
\tilde{i}_w\beta=\sum_{j=1}^q (-1)^{j-1}
i_{v_j}\beta\otimes 
v_1\wedge\ldots\wedge v_{j-1}\wedge\widehat{v_j}\wedge v_{j+1}
\wedge \ldots\wedge v_q   
\enspace ,\end{equation}
where as above $w=v_1\wedge v_2\wedge\ldots\wedge v_q$, together with
the assumption that $\tilde{i}$ is also $C^{\infty}(P)$-linear in $w$.
Extend $\pi^\#\colon \Omega^l(P)\otimes \CX^k(P)\to \CX^{l+k}(P)$
so that $\pi^\#(\beta\otimes w)=(\pi^\#\beta)\wedge w$. Then
(\ref{eq:Poisson-Lie-2}) turns into
\begin{equation}
\label{eq:Poisson-Lie}   
\CL_\alpha w   =L_\xi w+\pi^\#\tilde{i}_w\,d\alpha
\enspace ,\end{equation}
an equation similar to (\ref{eq:Poisson-Lie-1}).
In Lemma \ref{lemma:invariance} we shall use the following two simple 
properties of $\tilde{i}$:
\begin{equation}
\label{eq:ti-i}
\tilde{i}_w\alpha=i_\alpha w
\enspace , \end{equation}
where $\alpha\in \Omega^1(P)$ and $w\in \CX^\s(P)$, and
\begin{equation}
\label{eq:ti-wedge}
\tilde{i}_w(\alpha_1\wedge\alpha_2)=
(-1)^{(k-1)l}\alpha_2\wedge\tilde{i}_w\alpha_1
+(-1)^k\alpha_1\wedge\tilde{i}_w\alpha_2
\enspace , \end{equation}
where $\alpha_1\in\Omega^k(P)$ and $\alpha_2\in\Omega^l(P)$.

In the same vein, one can define the ``Lie derivative'' $\CL_\alpha{\beta}$
of a differential from $\beta\in \Omega^\s(P)$ in the direction
of a one-form $\alpha$. (See \cite{bv}, \cite{va:book}, and also
\cite{gi:moment} for a more detailed treatment.) 
When $\beta$ is a one-form, the derivative is given by (\ref{eq:bracket}):
$\CL_{\alpha}\beta=\{\alpha, \beta\}$. In general, $\CL_\alpha$ 
satisfies the Leibniz identity on forms (cf., (\ref{eq:Leibn})). Finally,
$\CL_\alpha$ of $w\in\CX^k(P)$ and $\beta\in\Omega^k(P)$ are 
related by the following formula from \cite{gi:moment}
\begin{equation}
\label{eq:dual}
\CL_{\alpha}\left<\beta,w\right>
=\left<\CL_{\alpha}\beta,w\right> + \left<\beta,\CL_{\alpha}w\right>
\enspace .
\end{equation}

Let now $P$ be acted on by a group $G$. Denote by $a\colon\fg\to\CX^1(P)$ 
the infinitesimal action homomorphism. The action is said to be 
{\em cotangential} 
if $a$ can be lifted to $\Omega^1(P)$ as a linear map, i.e., 
there exists a linear map $\tilde{a}\colon {\frak g}\to \Omega^1(P)$
such that $a=\pi^{\#}\circ \tilde{a}$. We call $\ta$
a {\em cotangent lift} or a {\em pre-momentum mapping}.
A pre-momentum mapping $\tilde{a}$ is said to be a {\em equivariant} if
$\tilde{a}$ is an anti-homomorphism of Lie algebras. 

It is clear that an action with an equivariant pre-momentum mapping
is cotangential and
a cotangential action is tangential, i.e., tangent to the symplectic
foliation of $\pi$. (The converse is not true \cite{gi:moment}.)
Note that here we do not require the action to
preserve $\pi$ or even to be Poisson. The problem of existence
and uniqueness of an (equivariant) pre-momentum mapping is analyzed 
in \cite{gi:moment} in detail. In particular, it is shown that 
{\em a cotangential action of a compact group with $H^2(\fg)=0$ admits an
equivariant pre-momentum mapping.}

Let $G$ be a Poisson Lie group and $G^*$ its dual group. 
Recall that $\pi_G$ vanishes at $e\in G$ and the linearization
of $\pi_G$ at $e$ is a linear homomorphism $\delta\colon \fg\to\fg\wedge\fg$,
which is a cocycle with respect to the adjoint action.
The Lie algebra of $G^*$ is $\fg^*$ with
the bracket defined as $\delta^*\colon\fg^*\wedge\fg^*\to\fg^*$ or
equivalently via (\ref{eq:bracket}). Thus every $\xi\in\fg$
can be thought of as an element of $(\fg^*)^*=T^*_e G^*$. Denote
by $\theta_\xi$ the left-invariant one-form on $G^*$ extending
$\xi$. Note that the dressing action of $\xi$ on $G^*$ is given
by the vector field $-(\pi_{G^*})^\#\theta_\xi$.

Recall that $G$ is said to act on $P$ is a Poisson way if the
action map $G\times P\to P$ is Poisson. The following is
a very convenient infinitesimal criterion for an action 
of a connected group to be Poisson \cite{lu-we:lie-poisson}:
{\em An action is Poisson if and only if, for all $\xi\in\fg$,}
\begin{equation}
\label{eq:Poisson-act}
d_\pi a(\xi)=a\delta(\xi)
\enspace .\end{equation}
Here we use the same notation $a$ for the infinitesimal action
$\fg\to\CX^1(P)$ and for the induced homomorphisms 
$\wedge^\s a\colon \wedge^\s\fg\to\CX^\s(P)$.

As defined by Lu in \cite{lu:moment}, a {\em momentum mapping} for a 
Poisson action is a map $\mu\colon P\to G^*$ such that
\begin{equation}
\label{eq:moment}
a(\xi)=-\pi^\#\mu^*\theta_\xi
\enspace .\end{equation}
A momentum mapping is said to be {\em equivariant} if it is such with
respect to the dressing action or, equivalently, if it is Poisson 
\cite{lu:moment}. 
For example, the identity map is an equivariant momentum mapping for 
the dressing action. An (equivariant) momentum mapping gives rise
to an (equivariant) cotangential lift
by the formula $\tilde{a}(\xi)=-\mu^*\theta_\xi$. Thus an action with
a momentum mapping is necessarily tangential and even cotangential.

The following theorem will play a crucial role in the subsequent analysis.
Recall that once a pre-momentum mapping $\tilde{a}$ is fixed, the spaces
$\CXP$ and $\Omega^\s(P)$ become $\fg$-modules via $\CL_{\tilde{a}(\cdot)}$. 
(Note that the structure of $\fg$-modules may depend
on the choice of $\tilde{a}$.)
Assume that $P$ is equipped with a genuine, not only infinitesimal, 
$G$-action.

\begin{thm}  
\label{thm:G-mod}   
~~~~~~~~
\begin{enumerate}
\item[{\rm (i)}]{\em \cite{gi:moment}}
The $\fg$-modules $\CXP$ and $\Omega^\s(P)$ integrate to representations the 
universal covering $\tilde{G}$ of $G$. The resulting $\tilde{G}$-action 
$\CXP$ commutes with $d_\pi$
and induces a trivial $\tilde{G}$-action on $\HP$ (see Section 
{\em \ref{sec:Poisson-coh}}). 

\item[{\rm (ii)}] Assume that $G$ is compact and $\ta$ is associated
with an equivariant momentum mapping. Then the infinitesimal action
integrates to representations on $\CXP$ and $\Omega^\s(P)$ of a finite 
covering $\bar{G}$ of $G$.

\end{enumerate}
\end{thm}

\proof Assertion (i) was proved in \cite{gi:moment}. 
By (\ref{eq:dual}), it is sufficient to prove (ii) for $\CXP$.
Recall that some finite covering $\bar{G}$ of $G$ can be decomposed as 
$\T\times K$, where $\T$ is a torus and $K$ is a compact simply 
connected (semisimple) Lie group. We claim that the infinitesimal action 
integrates to a representation of $\bar{G}$. 

A straightforward calculation with the cocycle $\delta$ shows that $\pi^G$
vanishes along $\T$, and so $\T$ is a Poisson subgroup of $G$. 
(Note that in contrast with $\T$, the subgroup $K$ does not, in general, 
inherit the structure of a Poisson Lie group, and we cannot apply (i) 
to it.) Furthermore, $\ft^*$ is an abelian subalgebra in $\fg^*$ and we 
have a Poisson projection
$G^*\to\ft^*$ (cf., \cite{lu-we:lie-poisson}). The composition of this
projection with the momentum mapping is an equivariant momentum mapping
for the induced $\T$-action. Clearly, $d\ta(\xi)=0$ for $\xi\in\ft$, 
which yields, by (\ref{eq:Poisson-Lie}), $L_{a(\xi)}=\CL_{\ta(\xi)}$ 
for $\xi\in\ft$. This shows
that the infinitesimal $\ft$-action on $\CXP$ and $\Omega^\s(P)$ integrates 
to a $\T$-action. Now it is a routine to conclude that the $\fg$-modules 
integrate to modules over $\T\times K$.\qed

We will also need the following result (cf., equation (\ref{eq:Poisson-act})).

\begin{lemma}  
\label{lemma:technical}
Let $\ta$ arise from an equivariant momentum
mapping. Then, for every $\xi\in\fg$,
\begin{equation}
\label{eq:d-delta}
d\ta(\xi)=\ta(\delta(\xi))
\enspace ,\end{equation}
where as before, the same notation $\ta$ is used for the pre-momentum
mapping and for the induced maps $\wedge^\s\fg\to\Omega^\s(P)$.
\end{lemma}

\proof It is clear that equation (\ref{eq:d-delta}) holds on $G^*$.
(See, e.g., \cite{lu:moment}.) On $P$, (\ref{eq:d-delta}) is just
the pull-back of the equation on $G^*$.\qed

\begin{rmk}
{\em
{\em  Assume that $G$ is compact or semisimple. Then
the cocycle $\delta$ is exact: there exists $r\in\wedge^2\fg$ such that
$\delta(\xi)=[\xi,r]$. As a consequence, $d\ta(\xi)=\CL_{\ta(\xi)}\ta(r)$,
when $\ta$ arises from an equivariant momentum mapping.}

It is well-known that $r$ exists when $\fg$ is semi-simple:
any one-cocycle on a finite-dimensional semisimple Lie algebra is
exact due to Whitehead's lemma. (See, e.g., \cite{gui:book}, Section II.11.) 
Assume that $G$ 
is compact. The function $\Delta\colon G\to\wedge^2\fg$,  defined as 
$\Delta(g)=(R_{g}^{-1})_*\pi_g$, where $R_g$ stands for the left 
translation, is a nonhomogeneous cocycle on $G$ with respect to the 
adjoint action. (See, e.g., \cite{lu-we:lie-poisson}
and references therein.) 
Since $G$ is compact, any cocycle is exact \cite{gui:book}. The cocycle 
$\delta$ is the image of $\Delta$ under the natural homomorphism from 
the cochains on $G$ to the cochains on $\fg$. Thus $\delta$ is exact 
as well. Furthermore, we can even produce an explicit 
expression for $r$ by applying a homotopy formula from \cite{gui:book}:
$
r=\int_G Ad_{g^{-1}}(R_{g}^{-1})_*\pi_g\,dg
$.
}
\end{rmk}

Rather little seems to be known on the existence of equivariant
momentum mappings when the Poisson structure on $G$ is nontrivial. 
(See a discussion in \cite{gi:moment}.) When $P$ is symplectic and
simply connected, nonequivariant momentum mappings $P\to G^*$
exist and are parameterized by elements of $G^*$ \cite{lu:moment}. 
If $G$ is semisimple, an equivariant momentum mapping $P\to G^*$
is at most unique \cite{gi:moment}.

\subsection{Poisson cohomology}
\label{sec:Poisson-coh}

In this section we recall the definition and some of the properties of 
the Poisson cohomology spaces, a notion introduced by Lichnerowitz in 
\cite{lich}. (A general introduction can be found, for example, in
\cite{va:book}. A detailed discussion more oriented toward our present 
goals is given in \cite{gi:moment}.)

Since $[\pi,\pi]=0$, the operation
$d_\pi=-[\cdot,\pi]\colon\CX^\s(P)\to\CX^{\s+1}(P)$ is a differential:
$d^2_\pi=0$. The {\em Poisson cohomology} $H^{\star}_{\pi}(P)$
is the cohomology of the complex $(\CX^{\star}, d_{\pi})$.
 
\begin{exam}
\label{exam:coh_sympl}
{\em
The homomorphism $\pi^\#\colon \Omega^\s(P)\to\CX^\s(P)$ commutes
with the differentials up to a sing, and so induces a homomorphism
$\pi^\#\colon H^\s(P)\to H^\s_\pi(P)$. When $P$ is symplectic,
$\pi^\#$ is an isomorphism on the level of complexes and in cohomology.
On the other hand, $\pi^\#=0$ when $\pi=0$, and 
$H^{\star}_{\pi}(P)=\CX^{\star}(P)$. In general, Poisson cohomology
shares the properties of de Rham cohomology and multi-vector fields
on $P$. In effect, the Poisson cohomology classes on $P$ can be thought
of as multi-vector vector fields on $P$ in the ``Poisson category'' 
\cite{gi-lu:calculus}.
}
\end{exam}

\begin{exam}
\label{exam:interpretaions-Poisson}
{\em
{\em Interpretations of Poisson cohomology:}

(0) $H_{\pi}^0$ is the algebra of the so-called Casimir functions,
i.e., functions constant on the leaves of the symplectic foliation.
Note that $H^{\star}_{\pi}(P)$ is an algebra over $H^0_{\pi}(P)$.

(1) $H^1_{\pi}(P)=Poiss/Ham$ is a Lie algebra (over $\reals$).
It is the quotient of the Lie algebra of $\pi$-preserving vector
fields, called Poisson, over the ideal of Hamiltonian
vector fields.

(2) $H^2_{\pi}(P)$ is the space of infinitesimal deformations
of $\pi$ modulo those given by infinitesimal diffeomorphisms.
}
\end{exam}

The following simple result is a very particular case of a general
interpretation due to Lu \cite{lu:hom} of the Poisson cohomology
of Poisson homogeneous spaces via certain Lie algebra cohomology.
When $\pi_G=0$, it was obtained by Koszul \cite{koszul}.

\begin{prop} 
\label{prop:Poisson-to-Lie} 
Let $G$ be  a Poisson Lie group and let $U$ be an open $G$-invariant
subset of $G^*$. Then
$$   
H^{\star}_{\pi}(U)=H^{\star}({\frak g}; C^{\infty}(U))   
\enspace ,$$   
where $C^{\infty}(U)$ is made into a ${\frak g}$-module by means 
of the dressing action of ${\frak g}$ on $G^{*}$.   
\end{prop}

In effect, we again have an isomorphism of complexes:
$\CX^{\s}(U)\simeq \wedge^{\s}\fg^*\otimes C^{\infty}(U)$.
When $G$ is compact, the cohomology can be calculated explicitly. 

\begin{thm}
\label{thm:Poisson-cohomology}
{\em \cite{gw}}.
Assume that $G$ is compact. Then 
$H^{\star}_{\pi}(U)=H^{\star}({\frak g})\otimes (C^{\infty}(U))^G$.
\end{thm}

The theorem follows immediately from Corollary \ref{cor:coh-functions} and
Proposition \ref{prop:Poisson-to-Lie} we just stated.

\section{Equivariant Poisson cohomology}
\label{sec:Poisson-cohomology}

\subsection{Basic properties of equivariant Poisson cohomology}

Let $(P,\pi)$ be a Poisson manifold with  an 
infinitesimal action of a Lie group $G$. Assume that the action admits 
an equivariant pre-momentum mapping $\tilde{a}\colon {\frak g}\to \Omega^1(P)$.

As we have shown above, $\CXP$ is a $\fg$-module with the $\fg$-action
defined via $\CL_{\ta(\cdot)}$.
It is easy to see that $\CXP$ equipped with the differential $d_\pi$
becomes a $\fg$-differential complex if we set 
$i_{\xi}w=i_{\tilde{a}(\xi)}w$ where $\xi\in\fg$ and $w\in\CXP$.
We call the $G$-equivariant cohomology of this complex the {\em 
equivariant Poisson  cohomology of $P$} and denote it by 
$H^{\star}_{\pi,G}(P)$. 

When $G$ is a Poisson Lie group and the action is Poisson and has an 
equivariant momentum mapping $\mu\colon P\to G^*$, 
the equivariant cohomology taken for $\tilde{a}(\xi)=-\mu^*\theta_\xi$
is said to be associated with $\mu$.

\begin{rmk}
{\em
We emphasize that $H^{\star}_{\pi,G}(P)$ is
defined regardless of whether we have a genuine or infinitesimal
$G$-action on $P$. It is also irrelevant for the definition, but not
for calculations, whether the action is Poisson or not. The only
data required is $(P,\pi)$, the action, and an equivariant pre-momentum
mapping $\ta$, which may exist even when the action is not
Poisson. For example, $\ta$ always exists and is unique when $P$ 
is symplectic. 
The complex $\CXP_{\fg}$ employed to define the equivariant cohomology 
depends on the choice of $\tilde{a}$ and $\ta$  is not unique in general. 
Furthermore, the very equivariant cohomology 
spaces appear to depend on $\tilde{a}$. (Thus the notation
$\CXP_{\tilde{a}}$ and $H^{\star}_{\pi,\tilde{a}}(P)$ would be more
appropriate.) However, it is not clear
whether this still can occur when $G$ is compact semisimple.

}
\end{rmk}

\begin{rmk}
{\em
Assume that $P$ is given a genuine $G$-action. By Theorem \ref{thm:G-mod}(i),
the structure of $\fg$-module on $\CXP$ via $\CL_{\tilde{a}}$ integrates to 
a $\tilde{G}$-module.

Moreover, it may happen that $\CXP$ is in effect a $G$-module.
This is the case, for example, when $P$ is symplectic or when the $G$-action 
preserves $\pi$ and $\tilde{a}$ arises from an equivariant
momentum mapping or just the one-forms $\ta(\xi)$ are closed. Knowing
that $\CXP$ is a $G$-differential complex may sometimes simplify the 
calculation of equivariant cohomology. 
When $G$ is compact and $\ta$ is associated with an 
equivariant momentum mapping, we can always assume from the
very beginning that $\CXP$ is a $G$-module due to Theorem 
\ref{thm:G-mod}(ii). Indeed, the replacement of $G$ by its
finite covering $\bar{G}$ has no effect on the equivariant cohomology.
}
\end{rmk}

\begin{exam}   
\label{exam:eq_coh_sympl}   
{\em   
As with ordinary Poisson cohomology, $\pi^{\#}$ induces a homomorphism
$$
H^{\s}_{G}(P)\longrightarrow H_{\pi,G}^{\star}(P)
\enspace .$$
Similarly to Example \ref{exam:coh_sympl}, 
 $\pi^{\#}\otimes id$ is an isomorphism of the complexes 
$\Omega^\s(P)_{\fg}\to\CX^\s(P)_{\fg}$ when $P$ is a symplectic,  and
$H_{\pi,G}^{\star}(P)\simeq H^{\star}_G(P)$.
}   
\end{exam}

The geometrical meaning of the equivariant Poisson cohomology
will be clarified in the next two sections. Here we just
point out that $H^\s_{\pi,G}(P)$ has nothing to do with
the Poisson cohomology of $P/G$ even when the action is free.
As shown in Section
\ref{subsec:tangent}, the equivariant Poisson cohomology is related
to the cohomology of the $G$-invariant multi-vector fields tangent to
the fibers of the momentum mapping.

\begin{exam}     
{\em    
Take $P=\complex^n\setminus 0$ with its standard  symplectic structure and 
the standard $G=\T^1$-action by rotations. Then 
$H^{\s}_{\pi,G}(P)\simeq H^\s_G(P)\simeq H^\s(P/G)= 
H^\s(\complex P^{n-1})$
because $P$ is symplectic. On the other hand, one can show that 
$H^{k}_{\pi}(P/G)=C^{\infty}(\reals^{+})$ when $k=0, ~2n-1$  
and zero otherwise.   
}   
\end{exam}

All the results of Section \ref{sect:poisson-general} apply to equivariant 
Poisson cohomology. In particular, $H^{\star}_{\pi,G}(P)$ is an algebra 
over $(S^\s\fg^*)^G$. 

\begin{exam}   
\label{exam:torus}  
{\em   
{\em Torus actions}. Let $G=\T^n$ be an $n$-dimensional torus   
$\T^1_1\times\ldots\times \T^1_n$ acting on $P$  with an  
equivariant pre-momentum mapping. Denote by  
$(\CXP)^{\tilde{a}}$ the subcomplex of $\CXP$ formed   
by multi-vector fields $w$ with $\CL_{\tilde{a}(\xi)}w=0$.  
Let us identify $S^{\s}\ft^*$ with $\reals[u_1,\ldots, u_n]$. As  
follows from the definition,  $H^{\s}_{\pi,\T}(P)$ is the cohomology of the 
complex   $\CXP^{\ft }\otimes S^\s\ft^*$ with the differential   
$d_{\pi}-\sum u_j i_{\tilde{a}(\xi_j)}$,  where $\xi_j$ is 
the generator of $\ft_j$.  
}  
\end{exam}

\begin{rmk}
{\em
In a similar way, one can introduce the algebroid equivariant cohomology.
The data required are an infinitesimal action of $G$ on the underlying
manifold $M$ of an algebroid $\CA$ and its lift to an action on $\CA$
satisfying the natural axioms. This is sufficient to make the 
complex $\wedge^\s\CA^*$ used in the definition of algebroid cohomology
into a $\fg$-differential algebra.
}
\end{rmk}

\begin{rmk}
{\em
It is interesting that $\fg$-differential complexes arise
in Poisson geometry in other ways as well. For example, as recently
pointed out by Jian-Hua Lu, the complex of Poisson homology $H^\pi_\s(P)$ with
the differential induced by the de Rham differential on $P$ is
a $\fg$-differential complex over the first Poisson cohomology
Lie algebra $\fg=H^1_\pi(P)$. More generally, $H^\pi_\s(P)$ is
a $\fg$-differential complex over the super Lie algebra $H^\s_\pi(P)$.
This is similar (and equivalent when $\pi=0$) to $\Omega^\s(P)$ being a 
$\fg$-differential complex over $\fg=\CX^1(P)$ or $\fg=\CX^\s(P)$.
Although the case of $\pi=0$ does not seem to yield any noteworthy results,
the general case may be more interesting.
}
\end{rmk}

\subsection{Equivariant Poisson cohomology of low degrees} 
\label{exam:interpretation}

Let us start by explicitly writing down the beginning of $\CX^{\star}_\fg$:
$$
0\longrightarrow
\CX^{0}(P)_\fg
\stackrel{d_G=d_{\pi}}{\longrightarrow}
\CX^{1}(P)_\fg\stackrel{d_G}{\longrightarrow}
\CX^{2}(P)_\fg\stackrel{d_G}{\longrightarrow}
\ldots
\enspace .$$
It is easy to see that $\CX^{0}(P)_\fg$ is the space of smooth 
$G$-invariant functions, since 
$L_{\pi^{\#}\tilde{a}(\xi)}=\CL_{\tilde{a}(\xi)}$ on $C^{\infty}(P)$.
Clearly $\Cas\subset C^{\infty}(P)^G$ and, in this 
degree, $d_G=d_{\pi}$. Thus
\begin{itemize}
\item{} $H^0_{\pi,G}(P)=H^0_{\pi}(P)=\Cas$ is just the algebra of Casimir
functions, and so $H^{\star}_{\pi,G}(P)$ is an algebra over $\Cas$.
\end{itemize}
Clearly, $\CX^1(P)_{\fg}$ is the space $\CX^1(P)^{\tilde{a}}$ 
of vector fields invariant with respect to $\CL_{\tilde{a}}$.
The differential on this space is
$d_Gw=d_{\pi}w+i_{\tilde{a}(\cdot)}w$. Finally, in the same
notation,  
$\CX^2(P)_{\fg}=\CX^2(P)^{\tilde{a}}\oplus (C^{\infty}(P)\otimes \fg^*)^G$.

\begin{thm} ~~~~
\label{thm:low-poisson}
\begin{enumerate} 

\item[{\em (i)}] Let $H^1(\fg)=0$. Then 
$$
H^1_{\pi,G}(P)
=\{w\in\CX^1(P)~\mid~d_\pi w=0~\hbox{\rm and}~ i_{\tilde{a}(\xi)}w=0
~\hbox{\rm for all}~\xi\in\fg\}/\,d_{\pi}(C^{\infty}(P)^G)
\enspace .$$

\item[{\em (ii)}] Assume that 
$\ta$ arises from an equivariant momentum
mapping $\mu\colon P\to G^*$. Then $H^1_{\pi,G}(P)$ is the quotient 
of the space of $G$-invariant Poisson vector fields tangent to the 
$\mu$-fibers over the space of Hamiltonian vector fields with 
$G$-invariant Hamiltonians.

\item[{\em (iii)}] Let $G$ be compact semisimple. Then the forgetful
homomorphism $H^\s_{\pi,G}(P)\to H^\s_{\pi}(P)$ is an isomorphism in degree
one and an epimorphism in degree two.
\end{enumerate}
\end{thm}

\begin{rmk} 
{\em
In other words, under the hypothesis of (ii), $H^1_{\pi,G}(P)$ is the quotient 
of vector fields whose (local) flows preserve $\pi$, the action, 
and $\mu$ over all such Hamiltonian vector fields.
Assertion (ii) gives an interpretation of 
in the spirit of Example \ref{exam:interpretaions-Poisson}.
}
\end{rmk}

\proof Assertion (i) is an immediate consequence of Theorem \ref{thm:low}.
To prove (iii), recall that by Theorem \ref{thm:G-mod}  the $\fg$-module 
structure on 
$\CX^\s(P)$ integrates to a $\bar{G}$-module for some finite covering
$\bar{G}$ of $G$. According to
our definition, $H^\s_{\pi,\bar{G}}(P)=H^\s_{\pi,G}(P)$.
The covering $\bar{G}$ is compact semisimple because so is $G$,
and (iii) follows from the second assertion of Theorem \ref{thm:low}.

Let us prove (ii). It suffices to show that in $\CX^1(P)^{\tilde{a}}$
every cohomology class can be represented by a $G$-invariant Poisson 
vector field. Since $\CL_{\tilde{a}(\xi)}w=0$, this is a consequence of 
the following observation:

\begin{lemma}
\label{lemma:invariance}
Assume that $\ta$ arises from an equivariant momentum mapping $\mu$.
Let $w\in\CX^\s(P)$ be tangent to the fibers of $\mu$, i.e.,
$i_{\tilde{a}(\xi)}w=0$ for all $\xi\in \fg$. Then
$\CL_{\tilde{a}(\xi)}w=L_{a(\xi)}w$ for all $\xi\in \fg$.
\end{lemma}

\proof First recall that
$$
\CL_{\tilde{a}(\xi)}w-L_{a(\xi)}w=\pi^\#\tilde{i}_w\,d\tilde{a}(\xi)
$$
due to (\ref{eq:Poisson-Lie}). Thus we just need to show that
the right hand side is zero.
Since $\tilde{a}$ arises from an equivariant momentum
mapping, Lemma \ref{lemma:technical} yields
$d\tilde{a}(\xi)=\tilde{a}(\delta(\xi))$
for $\delta(\xi)\in \wedge^2\fg$, and so
$$
\ta(\delta(\xi))=\ta(\xi_1)\wedge \ta(\zeta_1)+\ldots+
\ta(\xi_n)\wedge \ta(\zeta_n)
$$
for some elements $\xi_j$ and $\zeta_j$, $j=1,~2,\ldots, n$, of $\fg$.
According to equations (\ref{eq:ti-i}) and (\ref{eq:ti-wedge}),
$$
\tilde{i}_w\ta(\delta(\xi))=\sum \ta(\zeta_j)\otimes i_{\ta(\xi_j)}w
-\ta(\xi_j)\otimes i_{\ta(\zeta_j)}w
\enspace . $$
By the hypothesis, $i_{\ta(\xi_j)}w=i_{\ta(\zeta_j)}w=0$, and so
$\tilde{i}_w\ta(\rho_\xi)=0$, which completes the proof.\qed

\begin{rmk}  
{\em
It appears that the second equivariant Poisson cohomology  
should have an interpretation along the line  of Example 
\ref{exam:interpretaions-Poisson} as well.
Namely, one may hope that under certain conditions $H^{2}_{\pi,G}(P)$
becomes the tangent space to the moduli space of deformations of $\pi$,
the momentum mapping, and the action subject to some constraints. (By the 
moduli space we mean the quotient of the space of deformations by the action
of a group of diffeomorphisms.) At the moment, it is not clear to the 
author how to make such an interpretation rigorous.
}
\end{rmk}

Recall that when a Poisson Lie group $G$ acts on $P$ in a Poisson manner,
the graded space $\CX^\s(P)^G$ of $G$-invariant multi-vector fields on $P$
is closed under $d_\pi$ \cite{lu:hom}. (Here we take the standard $G$-action
on multi-vector fields. The differential $d_\pi$ need {\em not} commute 
with the action, but $\CX^\s(P)^G$ is still a subcomplex.)
The cohomology $H^\s(\CX^\s(P)^G,d_\pi)$ is called the
{\em invariant Poisson cohomology} of $P$. We will denote it by 
$H^\s_\pi(P)_G$. The inclusion of complexes induces a homomorphism
$j^\s\colon H^\s_\pi(P)_G\to H^\s_\pi(P)$.

\begin{exam}
\label{exam:coh-inv}
{\em  \cite{lu:hom}.
It is not hard to see that $H^\s_\pi(G^*)_{G^*}=H^\s(\fg)$. Here, as in
a number of similar examples considered before, we have an 
isomorphism of complexes: $\CX^\s(G^*)^{G^*}\simeq C^\s(\fg)$.
}
\end{exam}

\begin{cor}
Assume that $G$ is compact semisimple and the action of $G$ on $P$ admits
an equivariant momentum mapping. Then the
induced homomorphism $j^1\colon H^1_\pi(P)_G\to H^1_\pi(P)$ is surjective.
\end{cor}

\proof We need to show that every cohomology class in $H^1_\pi(P)$ can be 
represented by a $G$-invariant vector field. According to Theorem
\ref{thm:low-poisson}, it can be represented by a basic vector field $w$,
i.e., $w$ such that $d_\pi w=0$ and $i_{\ta(\xi)}w=0$ for all $\xi\in\fg$.
By Lemma \ref{lemma:invariance}, $w$ is $G$-invariant in the standard
sense. \qed

\begin{rmk}
{\em 
The corollary speaks in favor of a conjecture \cite{gi:moment}
that $j^\s$ is surjective when the group is compact and 
the action admits a momentum mapping.

Another conjecture on the homomorphism $j^\s$ having interesting 
applications is that $j^\s$ is injective (perhaps under some 
weak additional assumptions). When the action preserves $\pi$, both 
conjectures are proved in \cite{gi:moment}.
}
\end{rmk}

\subsection{Equivariant Poisson cohomology for locally free actions}
\label{subsec:tangent}

In this section we examine the equivariant Poisson cohomology
in the case where the momentum mapping is a submersion or,
more generally, the pre-momentum mapping is in a certain sense locally free.

Let $P$ be a Poisson manifold acted on in a Poisson fashion by a Poisson
Lie group $G$; also let, as usual, $\ta$ be an equivariant pre-momentum
mapping. Set
$$
\CX^\s(\ta)=\{w\in \CX^\s(P)\mid i_{\ta(\xi)}w=0~\hbox{\rm and}~
L_{a(\xi)}w=0~\hbox{\rm for all}~\xi\in\fg\}
\enspace .$$
When $\ta$ is associated with an equivariant moment map $\mu$, this
is just the space of all $G$-invariant multi-vector fields which
are tangent to the $\mu$-fibers. In this case, we will also use 
the notation $\CX^\s(\mu)$. Recall that by definition the basic
subcomplex of $\CX^\s(P)$ is 
$$
\CX^\s(P)_b=\{w\in \CX^\s(P)\mid i_{\ta(\xi)}w=0~\hbox{\rm and}~
\CL_{\ta(\xi)}w=0~\hbox{\rm for all}~\xi\in\fg\}
\enspace .$$
Our first objective is to compare the graded spaces $\CX^\s(\mu)$
and $\CX^\s(P)_b$.

\begin{prop}
\label{prop:basic-inv}
Assume that $\ta$ arises from an equivariant momentum mapping $\mu$. Then
$\CX^\s(\mu)=\CX^\s(P)_b $.
In particular, $\CX^\s(\mu)$ is a subcomplex of $\CX^\s(P)$ whose
cohomology denoted from now on by $H^\s_\pi(\mu)$, or $H^\s_\pi(\ta)$,
coincides with the basic cohomology $H^\s_{\pi,b}(P)$.
\end{prop}

The corollary is an immediate consequence of Lemma \ref{lemma:invariance}.

\begin{thm}
\label{thm:surj}
Assume that $G$ is compact and $\ta$ is associated with a 
momentum mapping $\mu$ which is a submersion onto its image. Then
$$
H^\s_{\pi,G}(P)=H^\s_{\pi,b}(P)=H^\s_\pi(\mu)
\enspace .$$
\end{thm}

According to the theorem, the Poisson equivariant cohomology 
has a very simple geometrical interpretation  when $\mu$ is surjective. 
Namely, it is just the cohomology of the complex
of multi-vector fields which are $G$-invariant and tangent to the
$\mu$-fibers.

Before we prove Theorem \ref{thm:surj}, let us state a more general 
result. An equivariant pre-momentum mapping $\ta$ is said to be 
{\em locally free} if it makes $\CX^\s(P)$ into a locally free 
$G$-differential algebra in the sense of Definition \ref{dfn:locally-free}.
More explicitly, this means that there exists a $\fg$-equivariant linear
map $\Theta\colon \fg^*\to\CX^1(P)$ such that 
\begin{equation}
\label{eq:equi}
\left<\ta(\xi),\Theta(\lambda)\right>=\lambda(\xi)
\end{equation} 
for all $\lambda\in\fg^*$ and $\xi\in\fg$. 

\begin{exam}
\label{exam:sym-free}
{\em
When $P$ is symplectic, $\ta$ is locally free if and only if
the action is locally free, i.e., the infinitesimal action
homomorphism $a\colon\fg\to\CX^1(P)$ is injective. Then 
$\pi^{\#}$ induces an isomorphism between the basic subcomplex  
$\Omega^\s(P)_b=\Omega^\s(P/G)$  of the de Rham complex 
and $\CX^\s(\mu)$. Hence, in particular,  
$H^\s(P/G)=H^\s_\pi(\mu)$.

In general,
it is easy to find an example where $\ta$ is locally free, but
the $G$-action on $P$ is trivial.
}
\end{exam}

As a particular case of Corollary \ref{thm:basic-algebra}, we have

\begin{thm}
\label{thm:poisson-loc-free}
Let $\ta$ be locally free. Then
$H^\s_{\pi,G}(P)=H^\s_{\pi,b}(P)$.
\end{thm}

\pproof of Theorem \ref{thm:surj}. By Theorem \ref{thm:poisson-loc-free}
it suffices to show that $\ta$ is locally free. 
For $\lambda\in\fg^*=T_eG^*$, denote by $\tilde{\lambda}$ the left-invariant
vector field on $G^*$ whose value at $e$ is $\lambda$. Then
$\theta_\xi(\tilde{\lambda})=\lambda(\xi)$.

Fix a connection on the fiber bundle $\mu\colon P\to G^*$
and let $\Theta(\lambda)$ be the horizontal lift of $\tilde{\lambda}$.
Then $
\left<\ta(\xi),\Theta(\lambda)\right>=\theta_\xi(\tilde{\lambda})
=\lambda(\xi)$.
Therefore, (\ref{eq:equi}) holds for $\Theta$. As mentioned
in Remark \ref{rmk:average-free}, we can make $\Theta$ be $\fg$-equivariant,
while keeping (\ref{eq:equi}), by averaging $\Theta$ over the 
$G$-action. \qed

\subsection{Calculations of equivariant Poisson cohomology}  
\label{subsec:calcul}

Let now $G$ be a Poisson Lie group and $K\subset G$ its closed
connected subgroup. By definition, the dressing action of $G$
on $G^*$ has an equivariant momentum mapping $\mu=id$. Thus the restriction
of the action to $K$ has a natural equivariant pre-momentum mapping, which
we denote by $\tilde{a}$ and use in the equivariant cohomology
$H^\s_{\pi,K}(U)$, where $U$ is a $G$-invariant subset
of $G^*$. The following result is an equivariant version 
of Proposition \ref{prop:Poisson-to-Lie}.

\begin{thm} 
\label{thm:equi-G-dual}
$H^\s_{\pi,K}(U)=H^\s(\fg,K;C^{\infty}(U))$.
\end{thm}

\begin{cor}
$H^0_{\pi,G}(U)=C^{\infty}(U)^G$ and $H^{\s>0}_{\pi,G}(U)=0$.
\end{cor}

\pproof of Theorem \ref{thm:equi-G-dual}. In the case at hands 
$\CX^\s(\mu)$ is just the Chevalley--Eilenberg complex $C^\s(\fg,K;C^{\infty})$.
Hence, the theorem is an immediate
consequence of Theorem \ref{thm:surj} which applies because $\mu=id$. 

Alternately, one may follow the line of the proof
of Proposition \ref{prop:Poisson-to-Lie}. The complex $(\CX^\s(U), d_\pi)$
can be naturally identified with the Chevalley--Eilenberg complex
$C^\s=(\wedge^\s\fg^*\otimes C^{\infty}(U), d_{Lie})$. Under
this identification, the Cartan model for the equivariant Poisson cohomology
of $U$ turns into the Cartan model for the $K$-equivariant cohomology
of the the Chevalley--Eilenberg complex $C^\s$ as in Section 
\ref{subsec:chev-eilen}. Hence, $H^\s_{\pi,K}(U)=H^\s_K(C^\s)$.
By Theorem \ref{thm:eq-Lie}, $H^\s_K(C^\s)=H^\s(\fg,K;C^{\infty}(U))$.
\qed

Combining Corollary \ref{cor:equi-expl} and Theorem \ref{thm:equi-G-dual},
we obtain

\begin{thm} 
Assume that $G$ is compact. Then
$H^\s_{\pi,K}(U)=H^\s(\fg,K)\otimes C^{\infty}(U)^G$.
\end{thm}

Finally, Theorem \ref{thm:spec-seq} implies the following

\begin{prop} 
Let $G$ be compact. Then
there exists a spectral sequence of $(S^\s\fg^*)^G$-modules 
with $E^{pq}_1=E^{pq}_2=H^q_{\pi}(P)\otimes (S^{p/2}\fg^*)^G$ which 
converges to $H^{\star}_{\pi,G}(P)$.

\end{prop}

\section {A spectral sequence associated with a momentum mapping}
\label{sec:spec-mom} 

Now we are ready to introduce and study the main object of
this paper -- a new spectral sequence 
associated with an equivariant momentum mapping. 
The examples are analyzed in Section \ref{sec:examples}.

\subsection{The spectral sequence}

Let $P$ be a Poisson manifold acted on in a Poisson fashion by a Poisson 
group $G$ and let the action admit an equivariant momentum mapping $\mu$.
The following result is an immediate consequence of Theorems 
\ref{thm:spec-new} and \ref{thm:G-mod}.

\begin{thm}
\label{thm:spec-poisson}
There exists a spectral sequence converging to $H^\s_\pi(P)$ and
having $E^{pq}_1=H^q(\fg)\otimes (\CX^\s(P)\otimes W^\s(\fg))^p_b$ and 
$E^{pq}_2=H^q(\fg)\otimes H^p_{\pi,G}(P)$.
\end{thm}

\begin{rmk}
{\em
For the sake of simplicity,  we deliberately chose to state the theorem
under rather restrictive assumptions.
The hypotheses of the theorem can be relaxed. Clearly, it is
sufficient to assume that we are given an equivariant pre-momentum
mapping for an action of some compact group $G$ (not necessarily Poisson)
on $P$ and that $\CXP$ integrates to a $G$-module.
}
\end{rmk}

As in Section \ref{sec:spec-general}, one can also consider a
decreasing filtration of $\CXP$ arising from the moment map:
\begin{equation}
\label{eq:filter-moment}
\CX^{p+q}_p(P)=\{ w\in \CX^{p+q}(P)\mid
i_{\alpha_1}i_{\alpha_2}\ldots i_{\alpha_{q+1}}w=0 
~\hbox{\rm for all}~ \alpha_1,\ldots,\alpha_{q+1}\}
\enspace ,
\end{equation}
where $\alpha_1,\ldots,\alpha_{q+1}$ are the pull-backs by $\mu$ 
of some one-forms on $G^*$. In effect, there exists a similar 
filtration for any Poisson map $P\to B$ of Poisson manifolds.
When $B=G^*$, it is sufficient to take only $\alpha_i$ which are
pull-backs of the left-invariant one-forms on $G^*$, i.e., 
$\alpha_i=\mu^*\theta_{\xi_i}$.
In general, there seems to be little chance to tell anything 
interesting about the spectral sequence for this filtration in addition 
to the fact that it converges 
to $H^\s_\pi(P)$. (See, however, Example \ref{exam:poiss5}.)
As in Theorem \ref{thm:spec-new-free},
the situation becomes much more favorable when $\CXP$ is a locally free 
$G$-differentiable algebra. As we know from the proof of Theorem
\ref{thm:poisson-loc-free}, this is the case when $\mu$ is a submersion
onto an open subset of $G^*$, which is then necessarily $G$-invariant.

\begin{thm}
\label{thm:spec-poisson-sub}
Assume that $G$ is compact and the momentum mapping $\mu$ is a 
submersion onto an open subset $U$ of $G^*$. The filtration
{\em (\ref{eq:filter-moment})} gives rise to a spectral sequence which
converges to $H^\s_\pi(P)$ and has $E^{pq}_1=H^q(\fg)\otimes\CX^p(\mu)$
and $E^{pq}_2=H^q(\fg)\otimes H^p_\pi(\mu)$ in the notation of Section
{\rm \ref{subsec:tangent}}.
\end{thm}

\begin{exam}   
\label{exam:spec-sympl}   
{\em   
Let $P$ symplectic. It is easy to see $\pi^\#$ induces an isomorphism    
from the spectral sequence of the principal $G$-bundle    
$M\times EG\to (M\times EG)/G$ (Example \ref{exam:spec-man})   
to that in Theorem \ref{thm:spec-poisson}.        

Assume also that $\mu$ is a submersion. Then according to Example   
\ref{exam:sym-free}, $\pi^{\#}$ induces an isomorphism between the basic   
de Rham complex $\Omega^\s(P)_b$ (e.g., $\Omega^\s(P/G)$    
if the action is free) and $\CX^\s(\mu)$. In particular, $\pi^{\#}$ gives   
rise to the identification $H^\s(P/G)=H^\s_\pi(\mu)$. Now we can add to   
it that $\pi^\#$ is an isomorphism between   
the spectral sequences of $M\to M/G$ (Example \ref{exam:spec-man-free})   
and the spectral sequence from Theorem \ref{thm:spec-poisson-sub}.   
We will elaborate on this example in the next section.           

Finally note that nontrivial systems of local coefficients do
not arise in these examples (nor in Theorems \ref{thm:spec-poisson} and
\ref{thm:spec-poisson-sub}) because $G$ is assumed to be connected.
}  
\end{exam}

\pproof of Theorem \ref{thm:spec-poisson-sub}. 
This result is an immediate consequence of Theorems
\ref{thm:spec-new-free} and \ref{thm:G-mod}, for the hypotheses guarantee
that $\CX^\s(P)$ is a locally free $G$-differential algebra. However,
since the proof of Theorem \ref{thm:spec-new-free} was omitted, we briefly
outline a direct proof of Theorem \ref{thm:spec-poisson-sub}
in order to illuminate the geometrical meaning and the nature 
of the initial terms and the differentials of the spectral sequence.
As mentioned above, the argument follows closely the proof of Theorem
\ref{thm:spec-new}. 

First recall that by the second assertion of Theorem \ref{thm:G-mod},
we may assume, replacing if necessary $G$ by its finite covering, that 
the $\fg$-action on $\CX^\s(P)$ via $\CL_{\ta}$ integrates to a $G$-action.

Denote, as above, by $\alpha_j$'s the pull-backs of some 
left-invariant forms on $G^*$. In other words, $\alpha_j=\mu^*\theta_{\xi_j}$
for some $\xi_j\in\fg$, and so the space of $\alpha_j$'s can be identified
with $\fg$. Observe that for $w\in\CX^{p+q}_p(P)$ 
the contraction
$$
c_w(\alpha_1,\ldots,\alpha_q)=i_{\alpha_q}\ldots i_{\alpha_{1}}w
$$
is a $p$-vector field tangential to the $\mu$-fibers, i.e.,
an element of $\CX^p_p(P)$. Hence, the correspondence $w\mapsto c_w$
gives rise to a linear homomorphism from $\CX^{p+q}_p(P)$ to the
space $C^q(\fg;\CX^p_p(P))$ of $q$-cochains on $\fg$ with coefficients
in $\CX^p_p(P)$. The kernel of this homomorphism is clearly 
$\CX^{p+q}_{p+1}(P)$. Thus $w\mapsto c_w$ descends to
a monomorphism, also denoted by $c$,
\begin{equation}
\label{eq:e_0}
E^{pq}_0=\CX^{p+q}_p(P)/\CX^{p+q}_{p+1}(P)\longrightarrow C^q(\fg;\CX^p_p(P))
\enspace .
\end{equation}
It is easy to see that $c$ is onto, and therefore an isomorphism.
(One can prove this using a distribution transversal to the 
$\mu$-fibers, which exists because  $\mu$ is a submersion.)
From now on we identify $E^{pq}_0$ and $C^q(\fg;\CX^p_p(P))$.

To find $d_0$ one uses an expression for $d_\pi$ which is analogous to
the Cartan formula for the de Rham differential but with the 
roles played by vector fields and differential forms interchanged. 
(See, e.g., \cite{bv} or formula (4.8) in Section 4.2 of \cite{va:book}.)
Then a straightforward calculation identical to that in the proof of
Theorem \ref{thm:spec-new} shows that $d_0$ is equal to the
Lie algebra differential $d_{Lie}$ on $C^\s(\fg;\CX^p_p(P))$.
Here we can take any one of two $\fg$-module structures, via $\CL$ or $L$, 
because by Lemma \ref{lemma:invariance} they coincide  on $\CX^p_p(P)$.
Thus we have
\begin{equation}
\label{eq:e_1}
E^{pq}_1=H^q(\fg;\CX^p_p(P))
\end{equation}
and, since $G$ is compact, $E^{pq}_1=H^q(\fg)\otimes (\CX^p_p(P))^G$
due to Theorem \ref{thm:Lie-cohomology}. By definition, 
$(\CX^p_p(P))^G=\CX^p(\mu)$, which yields
$$
E^{pq}_1=H^q(\fg)\otimes \CX^p(\mu)
\enspace .$$
(Recall that $\CX^p(\mu)$ is the space of $G$-invariant $p$-vector fields
which are tangent to the fibers of $\mu$; see Section \ref{subsec:tangent}.) 

It remains to show that $d_1$ coincides with the restriction of $d_\pi$
to $\CX^\s(\mu)$.  The rest of the proof is a routine, although
lengthy, analysis of $d_1$ serving to this end and based on the assumptions 
that $G$ is compact and $\mu$ is surjective. (Example \ref{exam:spec-sympl} 
lends itself as a hint.) 

Let $w\in F^{p+q}_p$ represent an element $[w]$ of 
$E^{pq}_1$. Note that $d_0w=0$ yields $d_\pi w\in F^{p+q+1}_{p+1}$, which in 
turn guarantees that 
$i_{\alpha_q}\ldots i_{\alpha_1} d_\pi w\in \CX^{p+1}(\mu)$.
By definition, $d_1[w]\in E^{p+1\: q}_1$ is represented by the cocycle 
$$
(\alpha_1, \ldots,\alpha_q)\mapsto i_{\alpha_q}\ldots i_{\alpha_1} d_\pi w
$$
with values in $\CX^{p+1}(\mu)$. To ensure that $d_1$ has the desired form,
it suffices to show that for some particular choice of $w$ in $[w]$, the same
up to a sign cocycle is obtained when the contractions with 
$\alpha_1\ldots \alpha_q$ are followed by $d_\pi$, i.e., 
$$
i_{\alpha_q}\ldots i_{\alpha_1} d_\pi w=
\pm d_\pi i_{\alpha_q}\ldots i_{\alpha_1} w
\enspace .$$
(Since $\CX^\s(\mu)$ is a subcomplex, the right hand side is indeed in
$\CX^{p+1}(\mu)$.) This amounts to proving that
\begin{equation}
\label{eq:d_1}
(i_{\alpha_q}\ldots i_{\alpha_1} d_\pi w)(\beta_1,\ldots,\beta_{p+1})=
(-1)^q(d_\pi i_{\alpha_q}\ldots i_{\alpha_1} w)(\beta_1,\ldots,\beta_{p+1})
\end{equation}
at every point $x\in P$, for any one-forms $\beta_1,\ldots,\beta_{p+1}$ on
the tangent space to the $\mu$-fiber through $x$. Thus let $\beta_j$ form
a basis in $T^*_x\mu^{-1}(\mu(x))$. To verify (\ref{eq:d_1}), we will
extend $\beta_j$'s to global one-forms near $x$ in a particularly convenient
way and then apply the Cartan formula for $d_\pi$. 

The extension is to be carried out so that to obtain the $G$-invariant 
one-forms with respect to the $G$-action
on $\Omega^1(P)$ arising from $\CL$. (See Theorem \ref{thm:G-mod}.)
Let us begin by extending $\beta_j$'s to $G$-invariant differential forms on
$\mu$-fibers defined in a neighborhood of $x$. To do so we first extend
$\beta_j$ along a slice transversal to the $G$-orbit through $x$ and then use
the translations by the $G$-action to extend the forms to a neighborhood
of $x$. (Note that along the $\mu$-fibers the derivatives $L$ and $\CL$ 
coincide and so do the two lifts of $G$-action to the tangent (and 
cotangent) bundles to the $\mu$-fibers.)
Finally, let us fix a $G$-invariant projection $K$ of $TP$ to the tangent
bundle to the $\mu$-fibers. Here the invariance is understood with
respect to the action arising from $\CL$. Such a projection exists
because the fibers have constant dimension and $G$ is compact. Taking
the composition of this projection with $\beta_j$'s, which have so far been
defined along $\mu$-fibers, we obtain genuine one-forms on a 
neighborhood of $x$. These forms, denoted again by $\beta_j$, are
$G$-invariant by (\ref{eq:dual}).

According to the representation of $E^{pq}_1$ as a tensor product, a
multi-vector field $w$ in the class $[w]$
can be chosen in the form $w=y_q\wedge v_p$, where
$v_p$ is a $p$-vector field tangent to the $\mu$-fibers and $y_q$ is a 
$q$-vector field which is the horizontal lift of a $G^*$-invariant
multi-vector field on $G^*$. The horizontality is understood, of course,
with respect to the ``connection'' $\ker K$.

Applying the Cartan formula for $d_\pi$ to evaluate the left hand side 
of (\ref{eq:d_1}), we obtain a sum of the terms involving the
derivatives $\CL_{\beta_j}$ and $\CL_{\alpha_i}$ and various pairwise
brackets of $\beta_j$ and $\alpha_i$. These terms can be divided into
three groups,
$$
(i_{\alpha_q}\ldots i_{\alpha_1} d_\pi w)(\beta_1,\ldots,\beta_{p+1})
=\hbox{\rm I}+\hbox{\rm II}+\hbox{\rm III}
\enspace , $$ 
as follows.

\begin{itemize}

\item The first group contains the terms with $\CL_{\beta_j}$ and the
brackets $\{\beta_l,\beta_k\}$. The sum for this group is exactly
the right hand side of (\ref{eq:d_1}):
\begin{eqnarray*}
\hbox{\rm I}&=&\sum_{1\leq l<k\leq p+1}(-1)^{2q+k+l-1}
w(\{\beta_l,\beta_k\},\alpha_1,\ldots,\alpha_q,\beta_1,\ldots,
\widehat{\beta_l},\ldots,\widehat{\beta_k},\ldots,\beta_{p+1})\\
&&\quad +
\sum_{1\leq j\leq p+1} (-1)^{q+j}
\CL_{\beta_j}w(\alpha_1,\ldots,\alpha_q,\beta_1,\ldots,
\widehat{\beta_j},\ldots,\beta_{p+1})\\
&=&(-1)^q(d_\pi i_{\alpha_q}\ldots i_{\alpha_1} w)(\beta_1,\ldots,\beta_{p+1})
\enspace . \end{eqnarray*}

\item The second group is made up of the terms with $\CL_{\alpha_i}$ and the
brackets $\{\alpha_l,\alpha_k\}$. Each of the summands in the group 
vanishes because $w$ is taken in the form $y_q\wedge v_p$:
\begin{eqnarray*}
\hbox{\rm II}&=&\sum_{1\leq l<k\leq q}(-1)^{k+l-1}
w(\{\alpha_l,\alpha_k\},\alpha_1,\ldots,
\widehat{\alpha_l},\ldots,\widehat{\alpha_k},\ldots,\alpha_q,
\beta_1,\ldots,\beta_{p+1})\\
&&\quad +
\sum_{1\leq j\leq p+1} (-1)^{j}
\CL_{\alpha_i}w(\alpha_1,\ldots,\widehat{\alpha_i},\ldots,\alpha_q,
\beta_1,\ldots,\beta_{p+1})\\
&=&0
\enspace . \end{eqnarray*}

\item The third group involves the terms with $\{\alpha_i,\beta_j\}$. 
The brackets vanish, for the forms $\beta_j$ are $G$-invariant,
and so $\{\alpha_i,\beta_j\}=\CL_{\alpha_i}\beta_j=0$:
$$
\hbox{\rm III}=\sum_{i,j}(-1)^{q+i+j-1}
w(\{\alpha_i,\beta_j\},\alpha_1,\ldots,
\widehat{\alpha_i},\ldots,\alpha_q, 
\beta_1,\ldots,\widehat{\beta_j},\ldots,\beta_{p+1})=0
\enspace .$$
\end{itemize}
This completes the proof of (\ref{eq:d_1}) and the proof of the
theorem. \qed

\begin{rmk}
{\em
The identifications (\ref{eq:e_0}) and (\ref{eq:e_1}) hold even when
$G$ is not compact, but $\mu$ is still a submersion. However, the
$E_2$-term is, in this case, much more difficult to calculate.

As clear from the proof, the cohomology $H^\s(\fg)$ can be understood 
with some extra insight as the invariant Poisson cohomology   
$H^\s_\pi(G^*)_{G^*}$ (Example \ref{exam:coh-inv}) or, when $P$ is   
symplectic, as $H^\s(G)$ (see Example \ref{exam:spec-sympl}).
}
\end{rmk}

\subsection{Examples and applications}
\label{sec:examples}

Throughout this section $G$ is assumed to be a compact connected
(Poisson) Lie group.

Assume first that $P$ is a symplectic manifold. Then,
as we noted in Example \ref{exam:spec-sympl}, 
$\pi^{\#}$ gives rise to an isomorphism between the
spectral sequence in the de Rham cohomology and that in the Poisson
cohomology making the latter particularly easy to write down.
We illustrate this point by two simple examples where the Poisson
structure on $G$ is assumed to be trivial.

\begin{exam}
\label{exam:sym1}
{\em
Let $G$ be the circle acting on the standard symplectic 
$\reals^{2n}=\complex^n$ by rotations
with Hamiltonian $\sum |z_j|^2/2$. Since the manifold is (equivariantly)
contractible, its (equivariant) cohomology is the same as of a point.
The spectral sequence from Theorem \ref{thm:spec-poisson}) has 
$E_2^{pq}=H^q(S^1)\otimes H^p(\complex P^\infty)$. It
converges to $H^\s(pt)$: the differential $d_2$ kills all the terms but 
$E^{00}_2=\reals$.

Take now $P=\complex^n\setminus 0$ with the same action. The momentum mapping
on $P$ is a submersion and we are in a position to apply Theorem
\ref{thm:spec-poisson-sub}. The spectral sequence converges to
$H^\s_\pi(P)=H^\s(S^{2n-1})$ and has 
$E^{pq}_2$ which is the tensor product of $H^q(S^1)$ and 
$H^\s_\pi(\mu)=H^\s(\complex P^{n-1})$. (The easiest way to see
the latter isomorphism is via $\pi^{\#}$.) Finally, $d_2$ kills
all the $E_2$-terms but $E^{00}_2=\reals$ and $E^{2n-1\: 1}_2=\reals$.
}
\end{exam}

\begin{exam}
\label{exam:sym2}
{\em
Let $P=T^*G$ with its natural symplectic 
structure and the (left) $G$-action. The spectral sequence converges to 
$H^\s_\pi(T^*G)=H^\s(G)$ by collapsing already in the $E_2$-term. We have
$E_2=H^\s(\fg)\otimes H^\s_\pi(\mu)$ with $H^\s(\fg)=H^\s(G)$ and
$H^\s_\pi(\mu)=H^\s(T^*G/G)=H^\s(pt)$.
}
\end{exam}

Let us now turn to the examples where $P$ is a genuine Poisson manifold.

\begin{exam}
\label{exam:poiss1}
{\em
Let $P$ be an open $G$-invariant subset of $G^*$ with 
$\mu=id\colon P\hookrightarrow G^*$ and the standard 
dressing action of $G$. Observe that
$H^q_\pi(\mu)=C^{\infty}(P)^G$ when $q=0$ and zero otherwise, since
$\mu=id$. Then the spectral sequence collapses in the $E_2$-term and
$H^\s_\pi(P)=H^\s(\fg)\otimes C^{\infty}(P)^G=E_2$ gives the decomposition
of $E_2$ as the tensor product. (Note that we thus obtain Theorem 
\ref{thm:Poisson-cohomology} as a corollary of Theorem 
\ref{thm:spec-poisson-sub}.)
}
\end{exam}

\begin{exam}
\label{exam:poiss2}
{\em
Consider $P=su(2)^*\setminus 0$. Let $\mu$ be a Casimir function
without critical points, e.g., the distance to the origin.
We have $P=S^2\times \reals_+$ with the Poisson structure
$\pi=t\pi_0$, where $\pi_0$ is the standard symplectic Poisson structure
on $S^2$ and $\mu=t$ is the coordinate on $\reals_+$. Clearly,
$\mu$ gives rise to the trivial circle action on $P$, and we can 
apply Theorem \ref{thm:spec-poisson-sub}. The resulting spectral sequence
converges to 
$H^\s_\pi(P)=H^\s(su(2))\otimes \CA=\CA\oplus [u\wedge\pi_0]\cdot\CA$,
where $u=\partial_t$ and for the sake of 
brevity the algebra of Casimir functions $C^\infty(\reals_+)$ is denoted 
by $\CA$. Since the level sets of $\mu$ are just the symplectic leaves,
$H^\s_\pi(\mu)$ can be identified with the {\em tangential}
Poisson cohomology $H^\s_{\pi,tan}(P)$ introduced in 
\cite{gi:moment}, i.e., the cohomology of multi-vector fields tangent 
to the symplectic leaves. It is not hard to see that 
$H^\s_{\pi,tan}(P)=
\CA\oplus [\pi_0]\cdot \CA$.
Then, fixing also the generators of $H^\s(S^1)$, we get
\begin{equation}
\label{eq:e_2}
E_2^{pq}\simeq \left\{\begin{array}{ll}
C^{\infty}(\reals_+) & \mbox{if $p=0,~2$ and $q=0,~1$}\\
0                    & \hbox{otherwise}
\end{array}
\right.
\enspace . 
\end{equation}
The differential $d_1$ is an isomorphism $E_2^{01}\to E_2^{20}$ canceling
these two terms and leaving $E^{00}_2$ and $E^{21}_2$ intact.

Of course, the Poisson cohomology in this and the next example has been
calculated repeatedly by many authors using various methods.
See, in particular, \cite{VK}, \cite{Xu} and references therein.
}
\end{exam}

\begin{exam}
\label{exam:poiss3}
{\em
In the notation of Example \ref{exam:poiss2}, let us take $P=S^2\times\reals_+$
with $\pi=f(t)\pi_0$, where $f$ is a nonvanishing function and
$\mu=t$ the projection onto $\reals_+$. (Of course, in what follows
$\reals_+$ can be replaced by $\reals$ or $S^1$.) We will
use the spectral sequence of Theorem \ref{thm:spec-poisson-sub} to
calculate the Poisson cohomology of $P$ and thus answer a question of
Gregg Zuckerman. It is easy to see that the $E_2$-term is independent
of $f$ and therefore given by (\ref{eq:e_2}).
As before, the spectral sequence stabilizes no later than in the 
$E_3$-term.
Moreover, $E^{00}_2$ and $E^{21}_2$ are not touched by $d_2$,
which reflects the fact that $H^0_\pi(P)$ and $H^3_\pi(P)$ are 
independent of $f$. The
only effect $f$ has is via $d_2\colon E_2^{01}\to E_2^{20}$. (For
example, $d_2=0$ when $f=const$ and, as we have seen, $d_2$ is an 
isomorphism when $f(t)=t$.)  One can show that under the identification
(\ref{eq:e_2}), $d_2\colon \CA \to \CA$ is
just the multiplication by $f'$. Thus denoting by $I$ the interior 
of the set $\{f'=0\}$, we get:
$$
H^k_\pi(P)=\left\{ \begin{array}{cl}
\CA=\Cas & \mbox{when $k=0$} \\
\left[ u \right] C^{\infty}_{0}(I) & \mbox{when $k=1$} \\
\left[\pi_0\right] \bigl(\CA/(f'\cdot \CA)\bigr) 
& \mbox{when $k=2$}\\
\left[ u\wedge\pi_0\right] \CA & \mbox{when $k=3$}
\end{array}
\right.
\enspace , $$
where $C^{\infty}_0(I)$ is the space of
functions on $I$ which, being set zero on the complement to $I$, extend 
smoothly to $\reals_+$. The requirement that $f$ is nowhere constant
renders $I=\emptyset$, and so $H^1_\pi(P)=0$. Furthermore, if $f$ is a 
Morse function with $n$ critical points,
$H^2_\pi(P)\simeq \CA/(f'\cdot \CA)=\reals^n$.
}
\end{exam}

\begin{exam}
\label{exam:poiss4}
{\em
In the setting of Example \ref{exam:poiss3}, let us take $P=S^2\times S^1$
so that $P$ become compact and $\CA=C^{\infty}(S^1)$. Also, take $f$ such that
$I$ is a proper subset of $S^1$, i.e., $I$ is neither $\emptyset$ nor $S^1$.
Then $C^\infty_0(I)$ is {\em not} a finitely generated $\CA$-module. As a
consequence, {\em $H^\s_\pi(P)$ is not a finitely generated module over
the algebra of Casimir functions even though $P$ is compact and the
symplectic foliation is just the direct product $S^2\times S^1$}. Of course,
$\pi$ itself fails to be real analytic even though its symplectic foliation
is such.
}
\end{exam}

Finally, let us analyze an example where we have only a Poisson
map but no group action.

\begin{exam}
\label{exam:poiss5}
{\em
Let $P=B\times X$ and let the symplectic leaves of $\pi$ be the
fibers $b\times X$, $b\in B$. The natural projection $\mu\colon P\to B$
is a Poisson submersion, provided that $B$ is equipped with the zero 
Poisson structure. As mentioned above, we still have the filtration
(\ref{eq:filter-moment}) associated with $\mu$ and thus a spectral
sequence converging to $H^\s_\pi(P)$. One may show that
$E^{pq}_2=\CX^q(B)\otimes H^p(X)$. This decomposition, unlike $d_2$,
is independent of $\pi$ as long as the symplectic foliation remains fixed.

It appears to be a safe and provable conjecture that $H^\s_\pi(P)$ is a 
finitely generated $C^\infty(B)$-module when $\pi$ is real analytic.

Assume that $B$ is an open subset of $\reals^n$. Then along the
line of Example \ref{exam:poiss3}, we may write 
$E^{pq}_2=H^q(\ft^n)\otimes H^\s_\pi(\mu)$, where $\mu$ is viewed as the
momentum mapping of a trivial $\T^n$-action. As before we have
$H^\s_\pi(\mu)=H^\s_{\pi,tan}(P)$, which can also be identified with
the tangential de Rham cohomology of the foliation into $\mu$-fibers.
}
\end{exam}

\section{Appendix: Proof of Theorem \ref{thm:Lie-cohomology}}

Throughout the proof we keep the notation of Section \ref{section:lie-alg}.
Thus $C^\s=C^\s(\fg,K;V)$.

Denote by $I\colon C^n\to C^n(\fg,K;V^G)$ the averaging over the $G$-action:
$$
(I\phi)(\xi_1,\ldots,\xi_n)=\int_G\rho(g)\phi(\xi_1,\ldots,\xi_n)\, dg
\enspace ,$$
where $\xi_1,\ldots,\xi_n$ are elements of $\fg$ and $\phi\in C^n$.
A direct calculation shows that $I$ commutes with the differential $d_\rho$,
i.e., $d_{\rho} I= I\, d_{\rho}$, even though, in general, the operators 
$\rho(g)$ do not. (This easily follows from the fact that 
$I(\rho_*(\xi)v)=0$ for any $\xi\in\fg$ and $v\in V$.)
To prove (i), it suffices to show that $I\phi-\phi$ is exact when $\phi$
is closed. 

Let us identify $C^{\s}$ with a certain subcomplex $\Omega^\s_\rho$ of
the de Rham complex 
$\Omega^{\s}$ of differential forms on $G/K$ with values in $V$.  
Denote by $T^*$ the action of $G$ on $\Omega^{\s}$ induced by the 
left translations. Consider the diagonal action $T^{\rho}$ of $G$ on 
$\Omega^{\s}=\Omega^{\s}(G/K)\hat{\otimes} V$, i.e., the action via $T^*$ on 
the first term and via $\rho$ on the second. 
Let us identify $C^\s$ with $\wedge^\s T^*_e(G/K)$ as vector paces. Then
the complex $C^{\s}$ is topologically
isomorphic to the subcomplex $\Omega^{\s}_{\rho}\subset \Omega^{\s}$ of 
$T^{\rho}$-invariant forms. The isomorphism 
$\Phi\colon C^{\s}\to \Omega^{\s}_{\rho}$ sends a cochain $\phi$ to
the form $\phi^{\rho}$ which is a unique extension
of $\phi$ to a $T^{\rho}$-invariant form.  
In particular, 
$\phi^{\rho}_{e}=\phi$ and the evaluation at $e$ is the inverse to 
$\Phi$. We do have an isomorphism
of complexes: (\ref{eq:dL}) coincides with the Cartan formula
for $d_{dR}$. From now on we identify $C^\s$ and $\Omega^\s_\rho$ 
and omit $\Phi$ from our notation. Hence $d_{dR}=d_\rho$ on $\Omega^\s_\rho$.

It is easy to see that the averaging operator 
$\omega\mapsto \int_G T^*_g\omega\, dg$ 
preserves $\Omega^\s_\rho$  and its
restriction to $\Omega^\s_\rho$ coincides with $I$. Thus we denote
this averaging $\Omega^\s\to \Omega^\s$ by $I$ again.
Similarly to the averaging of real-valued forms, $I$ commutes with $d_{dR}$   
(\cite{gui:book}, Appendix E). Furthermore, using Lemma E.1 of \cite{gui:book},
it is not hard to show that $I$ induces the identity 
map on the de Rham cohomology, i.e., $I(\omega)-\omega$ is exact 
when $\omega\in \Omega^\s$ is closed. 

Consider a continuous projection $P\colon \Omega^{\s}\to \Omega^{\s}_{\rho}$ 
defined by
$$
P(\omega)=\int_G T^{\rho}_g\omega\,dg
\enspace .$$
A direct though lengthy calculation shows $P$ is a homomorphism
of complexes, i.e., $P$ commutes with $d_{dR}$. The calculation is based
on the observation that on $\Omega^0$, i.e., on smooth functions, $P$
commutes with the Lie derivatives along the left-invariant vector fields.
(Note that, unlike $I$, the projection $P$ 
does not induce the identity homomorphism in the cohomology unless 
the action of $G$ on $V$ is trivial.) 

Now we are ready to complete the proof of the theorem.
To prove (i), pick $\phi\in Z^n\subset C^n$. Then
$I\phi-\phi=d_{dR}\beta$
for some form $\beta$. The left hand side of this equation is in 
the subcomplex $C^\s=\Omega^\s_\rho$ of $\Omega^\s$. Applying $P$, we get
$I\phi-\phi=d_{\rho} P\beta$, for $P\,d_{dR}=d_\rho P$.

To prove (ii) and (iii),
we need the following property of $\Omega^\s$ to get around the
standard usage of the Hodge theory, which cannot be applied to
$V$-valued forms. Namely, {\em there exists a continuous linear map 
$\CH_{dR}\colon\Omega^{n}\to\Omega^{n-1}$ such that 
$\CH_{dR}\mid_{B^{n}_{dR}}$ is a right-inverse to $d_{dR}$, 
i.e., $d_{dR}\CH_{dR}\mid_{B^{n}_{dR}}=id$}. This observation is
an immediate consequence of a theorem by Mostow that $\Omega^{\s}$ is a 
{\em strongly injective}. (See \cite{bw}, Proposition 5.4 of Chapter IX or
\cite{gui:book}, Section D.1.1 and Lemma E.1.) This means   
that for every $n$, the inclusion of $Z^n_{dR}=\ker d_{dR}$ into $\Omega^n$ 
and that of $\Omega^n/Z^n_{dR}$ into $\Omega^{n+1}$ admit continuous left    
inverses.

Let us set $\CH= P\CH_{dR}\mid_{C^\s}$. Being a composition
of continuous linear maps, $\CH$ is continuous and linear. It is now 
sufficient to show that $d_{\rho}\CH(\phi)=\phi$ when $\phi\in B^{n}$. 
This is equivalent to showing that $d_{dR}P\CH_{dR}(\omega)=\omega$ for any
$d_{dR}$-exact $\omega\in\Omega^{n}_{\rho}$. Since $P$ and $d_{dR}$ commute,
and $\CH_{dR}\mid_{B^n_{dR}}$ is a right-inverse $d_{dR}$, we have 
$d_{dR}P\CH_{dR}(\omega)=P\,d_{dR}\CH_{dR}(\omega)=P\omega=\omega$. 

Finally note that $\CH_{dR}$ is independent of $\rho$ 
and $P$ and $\Phi$ are both smooth in $\rho$. As a 
consequence, $\CH$ depends smoothly on $\rho$ which proves (iii). \qed

\begin{rmk}
\label{rmk:closed}
{\em
It is clear from the proof that $B^n$ is a closed  subspace of $C^n$ 
and hence a Fr\'{e}chet space. Indeed, 
$B^n=\Phi^{-1}(\Omega^n_{\rho}\cap B^n_{dR})$, where as shown above
$\Omega^n_{\rho}$ and $B^n_{dR}$ are closed in $\Omega^n$.
}
\end{rmk}


\begin{thebibliography}{ABC}

\bibitem[AB]{at-bo:moment_map}
Atiyah, M., Bott, R., The moment map and equivariant cohomology,
{\em Topology}, {\bf 23} (1984), 1-28.

\bibitem[BV]{bv}
Bhaskara, K.H., Viswanath, K., Calculus on Poisson manifolds,
{\em Bull. London Math. Soc.}, {\bf 20} (1988) 68-72.

\bibitem[BW]{bw}
Borel, A., Wallach, N., {\em Continuous Cohomology, Discrete
Subgroups and Representations of Reductive Groups},
Annals of Math. Studies, no. 94, Princeton University Press, 
Princeton, NJ, 1980.

\bibitem[Ca1]{cartan1}
Cartan, H., Notions d'alg\`{e}bre diff\'{e}rentielle; application aux
groupes de Lie et aux vari\'{e}t\'{e} o\`{u} op\`{e}re un group de Lie.
In {\em Colloque de Topologie}, C.B.R.M, Bruxelles, 1950, 15-27.

\bibitem[Ca2]{cartan2}
Cartan, H., La transgression dans up groupe de Lie et dans up
espace fibr\'{e} principal.
In {\em Colloque de Topologie}, C.B.R.M, Bruxelles, 1950, 57-71.

\bibitem[DKV]{DKV}
Duflo, M., Kumar, S., Vergne, M., {\em Sur la cohomologie \'{e}quivariante
des vari\'{e}t\'{e} diff\'{e}rentiables}, Asterisque, {\bf 205}, 1993.

\bibitem[Es]{es:application} 
Est, E.T. van, On algebraic cohomology concepts in Lie groups, I, II,
{\em Indag. Math.} {\bf 17} (1955), 225-233, 286-294.

\bibitem[Fu]{fu:cohomology}
Fuks, D.B., {\em Cohomology of infinite-dimensional Lie algebras},
Consultants Bureau, New York, 1986.

\bibitem[Gi1]{gi}
Ginzburg, V.L., Some remarks on symplectic actions of compact groups,
{\em Math. Z.} {\bf 210} (1992), 625-640.

\bibitem[Gi2]{gi:moment}
Ginzburg, V.L., Momentum mappings and Poisson cohomology, 
{\em Int. J. Math.}, {\bf 7} (1996), 329-358.

\bibitem[GL]{gi-lu:calculus}
Ginzburg, V.L., Lu, J.-H., Poisson cohomology of Morita 
equivalent Poisson manifolds, {\em IMRN}, {\bf 10} (1992), 
199-205.

\bibitem[GW]{gw}
Ginzburg, V.L., Weinstein, A., Lie-Poisson structure on some
Poisson Lie groups, {\em J. Amer. Math. Soc.}, {\bf 5} (1992),
445-453.

\bibitem[Gu]{gui:book}
Guichardet, A., {\em Cohomologie des groupes topologiques et des
alg\`{e}bres de Lie}, Cedic/Fernard Nathan, Paris, 1980.

\bibitem[GLS]{gls}
Guillemin, V., Lerman, E., Sternberg, S., {\em Symplectic Fibrations
and Multiplicity Diagrams}, to appear.

\bibitem[Ko]{koszul}
Koszul, J.L., Crochet de Schouten-Nijenhuis et cohomologie,
{\em Ast\'{e}risque, hors serie, Soc. Math. France}, Paris (1985), 257-271.

\bibitem[Li]{lich}
Lichnerovitz, A., Les vari\'{e}t\'{e}s de Poisson et leurs alg\`{e}bres
de Lie associ\'{e}es, {\em J. Diff. Geometry}, {\bf 12} (1977), 253-488.

\bibitem[Lu1]{lu:moment}
Lu, J.-H., Momentum mappings and reductions of Poisson actions, in
Dazord, P., Weinstein, A., (editors), {\em Symplectic geometry,
groupoids, and integrable systems}, S\'{e}minaire Sud Rhodanien
G\'{e}om\'{e}trie \`{a} Berkeley (1989), Springer-Verlag, New York, 1991.

\bibitem[Lu2]{lu:hom}
Lu, J.-H., Poisson homogeneous spaces and Lie algebroids associated to
Poisson actions, Preprint 1995.

\bibitem[LW]{lu-we:lie-poisson}
Lu, J.-H., Weinstein, A., Poisson Lie groups, dressing actions, and
the Bruhat decomposition, {\em J. Diff. Geometry}, {\bf 31} (1990), 501-526.

\bibitem[Va]{va:book}
Vaisman, I., {\em Lectures on the Geometry of Poisson Manifolds},
Progress in Math. 118, Birkhauser, Boston, 1994.

\bibitem[VK]{VK}
Vorob'ev, Yu.M, Karasev, M.V., Poisson manifolds and the Schouten
bracket, {\em Functional Analysis and its Applications},
{\bf 22} no. 1 (1988), 1-9.

\bibitem[Xu]{Xu}
Xu, P., Poisson cohomology of regular Poisson manifolds,
{\em Ann. Inst. Fourier, Grenoble}, {\bf 42} (1992), 967-988.
\end{thebibliography}
\end{document}